\documentclass[journal]{IEEEtran}

\usepackage[utf8]{inputenc}
\usepackage{color}
\usepackage{array}
\usepackage{verbatim}
\usepackage{float}
\usepackage{amsmath}
\usepackage{amsthm}
\usepackage{amssymb}
\usepackage{graphicx}
\usepackage{longtable}
\usepackage[unicode=true,
 bookmarks=false,
 breaklinks=false,pdfborder={0 0 1},colorlinks=false]
 {hyperref}
\hypersetup{
 colorlinks,bookmarksopen,bookmarksnumbered,citecolor=blue,urlcolor=blue}
\usepackage[]{apacite}
\makeatletter
 \let\oldforeign@language\foreign@language
 \DeclareRobustCommand{\foreign@language}[1]{%
   \lowercase{\oldforeign@language{#1}}}

 \let\oldforeign@language\foreign@language
 \DeclareRobustCommand{\foreign@language}[1]{%
   \lowercase{\oldforeign@language{#1}}}

\ifCLASSINFOpdf
\else
\fi

\newtheorem{defn}{Definition}

\newtheorem{lem}{Lemma}

\newtheorem{prop}{Proposition}
\newtheorem{proper}{Property}
\newtheorem{thm}{Theorem}
\newtheorem{rem}{Remark}

\newtheorem{assum}{Assumption}

\begin{document}
%
%
\onecolumn
\noindent\rule{18.1cm}{2pt}\\
\underline{To cite this article:}
{\bf{\textcolor{red}{Hashim A. Hashim, Sami El-Ferik, and Frank L. Lewis. "Adaptive synchronisation of unknown nonlinear networked systems with prescribed performance." International Journal of Systems Science 48, no. 4 (2017): 885-898.}}}\\
\noindent\rule{18.1cm}{2pt}\\

\noindent{\bf The published version (DOI) can be found at: \href{http://dx.doi.org/10.1080/00207721.2016.1226984}{10.1080/00207721.2016.1226984} }\\

\vspace{40pt}\noindent Please note that where the full-text provided is the Author Accepted Manuscript or Post-Print version this may differ from the final Published version. { \bf To cite this publication, please use the final published version.}\\

\textbf{
	\begin{center}
		Personal use of this material is permitted. Permission from the author(s) and/or copyright holder(s), must be obtained for all other uses, in any current or future media, including reprinting or republishing this material for advertising or promotional purposes.\vspace{60pt}\\
	\end{center}
\vspace{360pt}
}
\footnotesize{ \bf
	\vspace{20pt}\noindent Please contact us and provide details if you believe this document breaches copyrights. We will remove access to the work immediately and investigate your claim.
} 

\normalsize

\twocolumn
\title{Adaptive synchronization of unknown nonlinear networked systems with prescribed performance}

\author{Hashim~A.~Hashim$^*$, Sami~El-Ferik, and~Frank L. Lewis
\thanks{$^*$Corresponding author, H. A. Hashim is with the Department of Electrical and Computer Engineering,
University of Western Ontario, London, ON, Canada, N6A-5B9, e-mail: hmoham33@uwo.ca}
\thanks{S.~El-Ferik is with the Department of Systems Engineering, King Fahd University of Petroleum and Minerals, Dhahran, 31261, Saudi Arabia.}
\thanks{F. L. Lewis is with the Department of Electrical and Computer Engineering, UTA Research Institute, The University of Texas at Arlington 7300 Jack Newell Blvd. S, Ft. Worth, Texas 76118.}
}


\markboth{--,~Vol.~-, No.~-, \today}{Hashim \MakeLowercase{\textit{et al.}}: Nonlinear Stochastic Attitude Filter on the Special Orthogonal Group}
\markboth{}{Hashim \MakeLowercase{\textit{et al.}}: Nonlinear Stochastic Attitude Filter on the Special Orthogonal Group}

\maketitle

\begin{abstract}
 This paper proposes an adaptive tracking control with prescribed performance function for distributive cooperative control of highly nonlinear multi-agent systems. The use of such approach confines the tracking error  within a large predefined set to a predefined smaller set. The key idea is to transform the constrained system into unconstrained one through the transformation of the output error. Agents' dynamics are assumed unknown, and the controller is developed for a strongly connected structured network. The proposed controller allows all agents to follow the trajectory of the leader node, while satisfying the necessary dynamic requirements. The proposed approach guarantees uniform ultimate boundedness for the transformed error as well as a bounded adaptive estimate of the unknown parameters and dynamics. Simulations include two examples to validate the robustness and smoothness of the proposed controller against highly nonlinear heterogeneous multi-agent system with uncertain time-variant  parameters and external disturbances.
\end{abstract}

\begin{IEEEkeywords}
Prescribed performance, Transformed error, Multi-agents, Distributed adaptive control, Consensus, Transient, Steady-state error, Networked Systems, Distributed Adaptive Control, Robustness.
\end{IEEEkeywords}

\IEEEpeerreviewmaketitle{}

\section{Introduction}

%
%
%
%
\IEEEPARstart{I}{n} recent years, distributive cooperative control has gained popularity and attention among control researchers owing to its capabilities to mimic the social behavior of animals such as bees swarming, birds flocking, ants foraging, fish schooling and so forth. Indeed, the control scheme enables a group of agents to perform a task that can be daunting for an individual agent in a simpler and faster manner. Besides, the cooperation between agents allows information exchange, improving performance and increasing productivity through collaboration, which emulates the standard behavior in social groups. Cooperative control contributes to the betterments of many applications such as the control of autonomous mobile robot vehicles in energy and mineral explorations, space studies, surveillance and other areas. Agents are connected by a communication network and exchange useful information. In such case, they are considered as nodes. The group of agents may follow one or more real or virtual leaders. The network formed by all nodes creates a graph, which can be directed or undirected. Undirected graphs refer to no difference between nodes while directed graphs define the direction of the flow of information between the node and its neighbors.\\

In the literature, \cite{fax_information_2004} and \cite{ren_consensus_2005} can be considered among the first pioneering studies addressing the consensus in multi-agent systems. Several other scholars contributed to the development in this field. For instance, the consensus of passive nonlinear systems has been addressed in \cite{chopra_passivity-based_2006}. Many research work such as \cite{olfati-saber_consensus_2007}-\cite{zhou_distributed_2015} investigated node consensus of cooperative tracking problem. Distributed tracking control for linear heterogeneous agents of MIMO systems with parameter uncertainties was established in \cite{zhao_fully_2014}. Cooperative tracking control for a single node has been studied in \cite{ das_distributed_2010, cao_distributed_2012} and in the case of high-order dynamics in \cite{zhang_adaptive_2012}. The work in \cite{das_distributed_2010,zhang_adaptive_2012} developed a neuro-adaptive distributed control for heterogeneous agents with unknown nonlinear dynamics connected through a digraph. In \cite{das_distributed_2010}, the authors considered nodes with first-order dynamic. Later on, in \cite{zhang_adaptive_2012}, high order systems have been addressed. In all previous studies, the input function in the node dynamics was assumed known.\\

On the other hand, cooperative tracking control problems of systems with unknown input function have been studied \cite{theodoridis_direct_2012} and \cite{el-ferik_neuro-adaptive_2014}. In \cite{theodoridis_direct_2012}, neuro-adaptive fuzzy was proposed to approximate unknown nonlinear dynamics and input functions.  The centers of the output membership functions are determined based on off-line trials. A very fundamental assumption in all these studies is the one that considers the unknown nonlinear dynamics as well as the input function as linear in parameters (LIP) (see for instance \cite{lewis_neural_1998} or \cite{el-ferik_neuro-adaptive_2014}). The goal is to guarantee the ultimate stability of the tracking error. Recently, there are several studies that were published addressing different issues related to adaptive control of multi-agent systems. These challenges include actuator fault (see for instance \cite{na2013adaptive}, \cite{ na2014adaptive},\cite{tong2014fuzzy}, \cite{tong2015fuzzy}, \cite{li2015prescribed} and \cite{zhao2016neural}, switching network topology \cite{yang2016distributed}, Predictor-based adaptive control \cite{wang2016predictor}, etc. All these challenging practical issues could benefit from prescribed performance framework to guarantee performance. In particular, a practical implementation of neuro-adaptive prescribed performance control to compensate for friction using a turn table servo system has been reported in \cite{na2014adaptive}. More implementation of such control approaches is really needed.\\

The distributed control of multi-agents attempts to tackle unknown nonlinearities, unmodeled dynamics, uncertainties, and disturbances. Estimation of the closed loop characteristics such as transient and steady state error is almost impossible to be represented analytically \cite{bechlioulis2014robust2}. Alternatively, prescribed performance has been proposed as a means to seclude the error to an arbitrarily small set, where the convergence is constrained to a given range. The key idea in the approach is to transform the error from the restricted space to the unconstrained one. The following section gives necessary details about the method. At this stage, it is worth to mention that prescribed performance approach aims to satisfy the following objectives. The convergence error has to be less that the predefined value; the transformed error is bounded; the maximum overshoot is less than the prescribed constant; the system's controlled output is smooth; and the control signal is both bounded and smooth.\\

Developing a cooperative control approach for multi-agent systems with prescribed performance has many benefits. In this context, the specified performance ensures that the consensus output error starts within a large predefined set and then converges systematically into a predefined narrow set \cite{bechlioulis_robust_2008}. During transient and steady-state, the tracking error satisfies a known time-varying performance. Adaptive cooperative control with prescribed performance has then the ability to increase the robustness of the system's behavior and to reduce the control effort. The proper selection of the upper and lower bounds of the prescribed performance functions guarantees the convergence of error within predefined limits smoothly and systematically. In the literature, robust adaptive control with prescribed performance function for feedback linearizable systems has been designed in \cite{bechlioulis_robust_2008}. The design of neuro-adaptive controllers to handle unknown nonlinearities and disturbances has been considered in \cite{bechlioulis_adaptive_2009}-\cite{yang_adaptive_2015} for different applications. The application of prescribed performance scheme with neural approximation included strict-feedback systems \cite{bechlioulis_adaptive_2009}, affine systems \cite{wang_verifiable_2010}, high order nonlinear systems \cite{bechlioulis_low-complexity_2014}. Each of these studies considered different assumptions on the input matrix continuity. Further refinement of these results improved the neural network weights previously tuned using trial and errors to avoid neural nets in the controller design and redevelop adaptive control with prescribed performance based on fuzzy adaptive tuning in \cite{sun_fuzzy_2014} and model reference adaptive control in \cite{mohamed_improved_2014}. All previous studies considered a single autonomous system. Recently, \cite{bechlioulis2014robust2} designed a control of a platoon with unknown nonlinear dynamics. Thus, the agents are set in a straight line, and each node sees only the one in front. This represents a particular network structure and a special formation for the nodes.\\

In this work, we propose a robust adaptive distributive control with prescribed performance for a group of nodes connected through a directed communication graph with known topology. The control law is fully distributed based on the fact that the control law of each agent respects the graph's topology and uses only the allowed local neighborhood information. Thus, the leader does not communicate with all the nodes. In our work, we consider a general network form characterized by its $L$ and $ B$ matrices.  The formation of the multi-agent systems can be anything including platooning.  The synchronization error between nodes follows a prescribed performance to satisfy predefined characteristics imposed by the designer. Each node contains unknown nonlinear dynamics and time-varying uncertainties. The controller is developed to meet a predefined transient response and specific characteristics of the steady-state synchronization error for each node. The original form of the prescribed performance as originally proposed in \cite{bechlioulis_robust_2008} is modified in this work to overcome the chattering in the control signal due to the interaction between nodes caused by the consensus algorithm the exchanging state information between nodes. The new approach guarantees stable dynamics with non-oscillatory, limited and smooth control signal.\\

The rest of the paper is organized as follows. Section \ref{sec_2} presents preliminaries of graph theory. Problem formulation and the derivation of the local error synchronization equation follow in Section \ref{sec_3}. Section \ref{sec_4} contains the control law development as well as the stability proof of the connected graph. Simulations results are presented in Section \ref{sec_5}. Conclusion and future work are in Section \ref{sec_6}.\\

\textbf{Notations:}The following symbols are used throughout the paper.\\
\begin{tabular}{l l l}
	$|\cdot|$ &:& absolute value of a real  number;\\
	$\|\cdot\|$ &:& Euclidean norm of a vector;\\
	$\|\cdot \|_F$ &:&  Frobenius norm of a matrix; \\
	$tr\{\cdot\}$ &:&  trace of a matrix;\\
	$\sigma(\cdot)$ &:&  set of singular values of a matrix, with the \\
	& &  maximum singular value $\bar{\lambda}$ and the minimum\\
	& &   singular value $\underline{\lambda}$;\\
	$P > 0 $ &:& indicates that the matrix $P$ is positive definite;\\
	& &  $(P \geq 0)$ (positive semi-definite);\\
	$\mathcal{N}$ &:&  set $\{1, . . . , N\}$;\\
	$I_m$ &:& identity matrix of order $m$.\\
	${\bf \underline{1}}$ &:&  unity vector $[1, \ldots, 1]^{\top} \ \in \mathbb{R}^n $ where $n$ is the\\
	& &  required appropriate dimension.
\end{tabular}

\section{Basic graph theory}\label{sec_2}
A graph is denoted by $\mathcal{G} = (\mathcal{V}, \mathcal{E})$ with a nonempty finite set of nodes (or vertices) $\mathcal{V} = \{\mathcal{V}_1, \mathcal{V}_2, \ldots, \mathcal{V}_n\}$, and a set of edges (or arcs) $E\subseteq \mathcal{V}\times \mathcal{V}$. $(\mathcal{V}_i,\mathcal{V}_j) \in E$ if there is an edge from node $i$ to node $j$. The topology of a weighted graph is often described by the adjacency matrix $A=[a_{ij}]\in \mathbb{R}^{N\times N}$  with weights $a_{ij} > 0$ if $(\mathcal{V}_{j}, \mathcal{V}_{i}) \in E$: otherwise $a_{ij} = 0$. Throughout the paper, the topology  of the communication network is fixed, i.e. $A$ is time-invariant, and the self-connectivity element $a_{ii} = 0$. A graph can be directed or undirected. A directed graph is called digraph. The weight in-degree of a node $i$ is defined as the sum of i-th row of $A$, i.e., $d_i=\sum_{j=1}^Na_{ij}$. Define the diagonal in-degree matrix $D = diag(d_1, \ldots, d_N)\in \mathbb{R}^{N\times N}$ and the graph Laplacian matrix $L = D-A$. The set of neighbors of a node $i$ is $N_i = \{j|(\mathcal{V}_j\times \mathcal{V}_i)\in E\}$. If node $j$ is a neighbor of node $i$, then node $i$ can get the information from node $j$ , but not necessarily vice-versa. For an undirected graph, the neighborhood requires a mutual relation. A direct path from node $i$ to node $j$ is a sequence of successive edges in the form $\{(\mathcal{V}_i,\mathcal{V}_k),(\mathcal{V}_k,\mathcal{V}_l), \ldots , (\mathcal{V}_m,\mathcal{V}_j\}$. A digraph has a spanning tree if there is a node (called the root) having a possible direct path to every other node in the graph. A digraph is strongly connected if for any ordered pair of nodes $[\mathcal{V}_i,\mathcal{V}_j]$ with $i\neq j$, there is a directed path from node $i$ to node $j$ (for more details, see \cite{ren_distributed_2008} or \cite{lewis_cooperative_2013}).

\section{Problem Formulation}\label{sec_3}
Consider the following nonlinear dynamics for the $i$th node
\begin{equation}
\label{eq:eq1}
\begin{aligned}
& \dot{x}_i = A_{mi}x_{i} + B_{mi}u_i + f_i\left(x_{i}\right) + w_i
\end{aligned}
\end{equation}
where the state node is $x_{i} \in \mathbb{R}$, the control signal node is $u_i \in \mathbb{R}$ and the unknown disturbance for each node is $w_i \in \mathbb{R}$. $A_{mi} \in \mathbb{R}$ and $B_{mi} \in \mathbb{R}$ are known constants; $f_i\left(x_{i}\right) \in \mathbb{R}$ is the unknown nonlinear part of the dynamics and assumed to be Lipschitz. From \eqref{eq:eq1}, the global dynamic can be written as
\begin{equation}
\label{eq:eq2}
\begin{aligned}
& \dot{x} = A_m\,x + B_m\,u  + f\left(x\right) + w
\end{aligned}
\end{equation}
where $x = [x_1,\ldots,x_N]^{\top}\in \mathbb{R}^{N}$, $u = [u_1,\ldots,u_N]^{\top}\in \mathbb{R}^{N}$, $f\left(x\right) = [f_1(x_1),\ldots,f_N(x_N)]^{\top}\in \mathbb{R}^{N}$, $w = [w_1,\ldots,w_N]^{\top}\in \mathbb{R}^{N}$, $A_m = {\rm diag}\{A_{m1},\ldots,A_{mN}\}$, and $B_m = {\rm diag}\{B_{m1},\ldots,B_{mN}\}$. $x_{0}$ is the leader's state and it represents the desired synchronization trajectory according to the following equation
\begin{equation}
\label{eq:eq3}
\begin{aligned}
& \dot{x}_0 = A_m\,x_{0}\,+\,f(x_{0},t)
\end{aligned}
\end{equation}
where $x_{0} \in \mathbb{R}$ is the leader state node, $f(x_{0},t) \in \mathbb{R}$ is the nonlinear part of the leader node's dynamic. The model presented in (\ref{eq:eq2})-(\ref{eq:eq3}) is very similar to the one treated in \cite{das_distributed_2010}, with the exception that the dynamic of the agent as described in (\ref{eq:eq2}) is more representative of real systems.
The local synchronization error function for agent $i$ can be described as in \cite{li_pinning_2004}-\cite{khoo_robust_2009}.
\begin{equation}
\label{eq:eq4}
\begin{aligned}
& e_i = \sum_{j \in N_i} a_{ij}(x_{i}-x_j)+b_{i}(x_{i}-x_{0})
\end{aligned}
\end{equation}
where $a_{ij} \geq 0$ and $a_{ii} = 0$. $a_{ij} > 0$ if agent $i$ is directed to agent $j$, $b_{i}\geq 0$. The network is such that $b_{i}> 0$ for at least one agent $i$. Hence, equation \eqref{eq:eq4} can be written in the global error form as
\begin{equation}
\label{eq:eq5}
\begin{aligned}
e &  = -\left(L+B\right)(\underline{1}x_{0}-x) = \left(L+B\right)(x - \underline{1}x_{0})\\
&  = \left(L+B\right)(\tilde{x})
\end{aligned}
\end{equation}
where the global error is $e = [e_1,\ldots,e_N]^{\top}\in \mathbb{R}^{N}$, global state vector is $\underline{1}x_{0} = \underline{x}_0 \in \mathbb{R}^{N}$, the Laplacian matrix is $L \in \mathbb{R}^{N \times N}$, $B \in \mathbb{R}^{N \times N}$ with $B = {\rm diag}\{b_i\}$ and $\underline{1}=[1,\ldots,1]^{\top} \in \mathbb{R}^{N}$. Note that $\tilde{x} = x - \underline{1}\cdot x_{0}$, and $\underline{f}(x_{0},t)=\underline{1}\cdot f(x_{0},t)$. For more details, the proof of equation \eqref{eq:eq5} is stated in \cite{lewis_cooperative_2013}.\\
The derivative error dynamics of \eqref{eq:eq5} is
\begin{equation}
\label{eq:eq7}
\begin{aligned}
& \dot{e} = \left(L+B\right)(A_m\, \tilde{x}+f\left(x\right)+B_m\,u+w-\underline{f}(x_{0},t))
\end{aligned}
\end{equation}
\begin{rem}
	The communication graph is considered strongly connected. Thus, if $b_i \neq 0$ for at least one $i$, $i=1,\ldots, N$ then $\left(L+B\right)$ is an irreducible diagonally dominant M-matrix and hence nonsingular \cite{Qu2009}.
\end{rem}

\begin{rem} (see \cite{lewis_cooperative_2013})
	If agent state is $x_{i} \in \mathbb{R}^{n}$ and the leader state $x_{0} \in \mathbb{R}^{n}$ where $n > 0$, then $e,x \in \mathbb{R}^{nN}$ and equation \eqref{eq:eq5} will be
	\begin{equation}
	\label{eq:eq6}
	\begin{aligned}
	& e = \left(\left(L+B\right)\otimes \mathbb{I}_N\right)(x - \underline{1}x_{0})
	\end{aligned}
	\end{equation}
	where $\otimes$ is the Kronecker product.
\end{rem}

Also, one should note that $B \neq 0$ for a strongly connected graph with
\begin{equation}
\label{eq:eq8}
\begin{aligned}
& ||\beta|| \leq ||e||/\underline{\lambda}\left(L+B\right)
\end{aligned}
\end{equation}
where $\underline{\lambda}\left(L+B\right)$ is the minimum singular value of $\left(L+B\right)$ \cite{lewis_cooperative_2013}.\\

A performance function $\rho\left(t\right)$ is associated with the error component $e\left(t\right)$ and is defined as a smooth function such as $\rho\left(t\right): \mathbb{R}_{+} \to \mathbb{R}_{+}$ is a positive decreasing function $\lim\limits_{t \to \infty}\rho\left(t\right)=\rho_{\infty}>0$. The prescribed performance function can be written as
\begin{equation}
\label{eq:eq9}
\rho_i\left(t\right)=(\rho_{i0} - \rho_{i\infty})\exp^{-\ell_i\, t}+\rho_{i\infty}
\end{equation}

where $\rho_{i0},\rho_{i\infty}$ and $\ell_i$ are appropriately defined positive constants. In order to overcome the difficulty caused through the synchronization algorithm and achieve the desired prescribed performance, the following time varying constraints are proposed:
\begin{equation}
\label{eq:eq10b}
-\delta_i\rho_i\left(t\right)<e_i\left(t\right)<\rho_i\left(t\right),\hspace{10pt} if \: e_i\left(t\right)>0
\end{equation}
\begin{equation}
\label{eq:eq11b}
-\rho_i\left(t\right)<e_i\left(t\right)<\delta_i\rho_i\left(t\right),\hspace{10pt} if \: e_i\left(t\right)<0
\end{equation}
for all $ t \geq 0 $ and $ 0 < \delta_i \leq 1 $, and $i=1, ...,N$.

\begin{rem}
	The dynamic constraints \eqref{eq:eq10b} and \eqref{eq:eq11b} represent a modification of the ones in \cite{bechlioulis_robust_2008}, and \cite{mohamed_improved_2014}. In these papers, the constraints are conditioned on $e(0)$ as follows
	\begin{equation}
	\label{eq:eq10}
	-\delta \rho \left(t\right)<e\left(t\right)<\rho \left(t\right),\hspace{10pt} if \: e(0)>0
	\end{equation}
	\begin{equation}
	\label{eq:eq11}
	-\rho \left(t\right)<e \left(t\right)<\delta \rho \left(t\right),\hspace{10pt} if \: e(0)<0
	\end{equation}
	Figure~\eqref{fig:fig1} shows the tracking error of controller with prescribed performance as it transits from a large to a smaller set in accordance with equations \eqref{eq:eq10b} and \eqref{eq:eq11b}.
	\begin{figure*}[ht]
		\centering
		\includegraphics[scale=0.4]{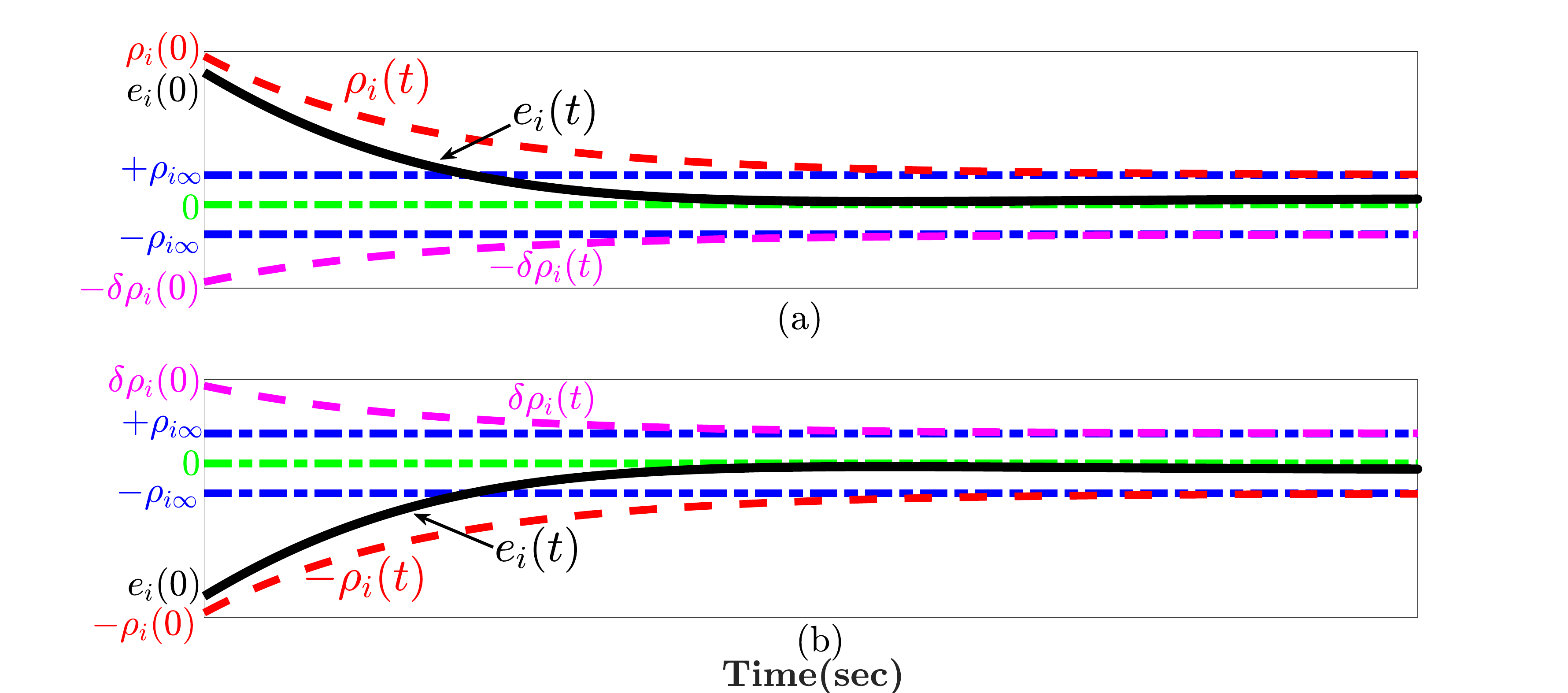}
		\caption{ Graphical representation of tracking error with prescribed performance
			(a) Prescribed performance of \eqref{eq:eq10}; (b) Prescribed performance of \eqref{eq:eq11}.}
		\label{fig:fig1}
	\end{figure*}
\end{rem}
Due to the interaction between agents' dynamics, such constraints will lead to instability. Upon crossing this reference, the system becomes unstable under the original formulation \eqref{eq:eq10} and \eqref{eq:eq11}. However, the switching based on $e_i\left(t\right)$ provides the necessary control to keep the system stable.
\begin{rem}
	In the control with prescribed performance framework as presented as in \cite{bechlioulis_robust_2008, mohamed_improved_2014}, the knowledge of the sign of $e_i(0)$ is sufficient to maintain the same robust controller  for all $t>0$ and satisfy the performance constraints (no switching occurs after $t=0$).  However, in the case of multi-agent systems, the synchronization error (\ref{eq:eq6}) creates a coupling between the different states of each agent. The interactions created may force the synchronization error to violate the desired performance constraints and exit the compact set if keep the same controller based solely on the sign $e_i(0)$, $i=1,\ldots,N$. Switching is rather needed at any time $t$ to keep the error within the compact sets.
\end{rem}

One should notice from Fig. \eqref{fig:fig1} that the tracking error in the case of multi-agent systems may exceed the lower (or upper) bound (in green color). Upon crossing this constraint, the system becomes unstable under the control based on $e_i(0)$. However, the switching based on $e_i\left(t\right)$ provides the necessary control effort to keep the system stable.

In order to transform the constrained error of the nonlinear system \eqref{eq:eq10b} and \eqref{eq:eq11b} to an unconstrained one, a transformed error $\epsilon_i$ is defined as
\begin{equation}
\label{eq:eq12b}
\epsilon_i=\psi\left(\frac{e_i\left(t\right)}{\rho_i\left(t\right)}\right)
\end{equation}
or equivalently,
\begin{equation}
\label{eq:eq131}
e_i\left(t\right)=\rho_i\left(t\right)S(\epsilon_i)
\end{equation}
where $\epsilon_i$, $S_i(.)$ and $\psi_i^{-1}(.)$ are all smooth functions, $i=1,2,\ldots,N$. $S(.)=\psi^{-1}(.)$ and $S(.)$ satisfy the following properties:
\begin{proper}
	\label{propertiesS}
	\begin{enumerate}
		\item $S_i(\epsilon_i)$ is smooth and strictly increasing.
		\item $-\underline{\delta}_i<S(\epsilon_i)<\bar{\delta}_i,\hspace{10pt} if \: e_i\left(t\right) \geq 0$\\
		$-\bar{\delta}_i<S(\epsilon_i)<\underline{\delta}_i,\hspace{10pt} if \: e_i\left(t\right)<0$
		\item ~\\
		$ \left.
		\begin{aligned}
		lim_{\epsilon_i \rightarrow -\infty}S(\epsilon_i)=-\underline{\delta}_i\\
		lim_{\epsilon_i \rightarrow +\infty}S(\epsilon_i)=\bar{\delta}_i,
		\end{aligned}
		\right\}
		\qquad {\rm if} \hspace{10pt} e_i\left(t\right)\geq 0$\\
		$\left.
		\begin{aligned}
		lim_{\epsilon_i \rightarrow -\infty}S(\epsilon_i)=-\bar{\delta}_i\\
		lim_{\epsilon_i \rightarrow +\infty}S(\epsilon_i)=\underline{\delta}_i,
		\end{aligned}
		\right\}
		\qquad {\rm if} \hspace{10pt} e_i\left(t\right)< 0$
	\end{enumerate}
\end{proper}
where\\
\begin{equation}
\label{eq:eq14b}
S(\epsilon_i)=
\left\{
\begin{aligned}
\frac{\bar{\delta}_i\exp^{\epsilon_i}-\underline{\delta}_i\exp^{-\epsilon_i}}{\exp^{\epsilon_i}+\exp^{-\epsilon_i}},& \quad \bar{\delta}_i> \underline{\delta}_i \hspace{5pt} {\rm if} \hspace{5pt} e_i\left(t\right)\geq 0\\
\frac{\underline{\delta}_i\exp^{\epsilon_i}-\bar{\delta}_i\exp^{-\epsilon_i}}{\exp^{\epsilon_i}+\exp^{-\epsilon_i}},& \quad \underline{\delta}_i > \bar{\delta}_i \hspace{5pt} {\rm if} \hspace{5pt} e_i\left(t\right) < 0\\
\end{aligned}
\right.
\end{equation}
Now, consider the general form of the smooth function
\begin{equation}
\label{eq:eq15b}
S(\epsilon_i)=
\begin{aligned}
& \frac{\bar{\delta}_i\exp^{\epsilon_i}-\underline{\delta}_i\exp^{-\epsilon_i}}{\exp^{\epsilon_i}+\exp^{-\epsilon_i}}
\end{aligned}
\end{equation}
and the transformed error
{\small
\begin{equation}
\label{eq:eq16b}
\begin{aligned}
\epsilon_i = & S^{-1}\left(\frac{e_i\left(t\right)}{\rho_i\left(t\right)}\right)\\
= & \frac{1}{2}\left\{ \begin{aligned} & \ln \frac{\underline{\delta}_i+e_i\left(t\right)/\rho_i\left(t\right)}{\bar{\delta}_i-e_i\left(t\right)/\rho_i\left(t\right)}, & {\rm with}\;\; \bar{\delta}_i > \underline{\delta}_i \hspace{5pt} {\rm if} \hspace{5pt} e_i\left(t\right)\geq 0 \\
& \ln \frac{\underline{\delta}_i+e_i\left(t\right)/\rho_i\left(t\right)}{\bar{\delta}_i-e_i\left(t\right)/\rho_i\left(t\right)} ,& {\rm with}\;\;\quad \underline{\delta}_i > \bar{\delta}_i \hspace{5pt} {\rm if} \hspace{5pt} e_i\left(t\right)< 0\\ \end{aligned}  \right.
\end{aligned}
\end{equation}
}
\begin{rem}
	In the previous set of equations, $\underline{\delta}_i$ and $\bar{\delta}_i$ exchange values depending on the sign of $e_i\left(t\right)$. One should note that the highest value of both involves subtracting the absolute value of $e_i\left(t\right)/\rho_i\left(t\right)$ and the lowest includes the addition of the absolute value of $e_i\left(t\right)/\rho_i\left(t\right)$.
	\label{remark5}
\end{rem}
Using this remark and the fact that $\rho_i\left(t\right)>0$, equation (\ref{eq:eq16b}) can be rewritten as
\begin{equation}
\label{eq:eq56b}
\begin{aligned}
\epsilon_i =
& \frac{1}{2}\left\{ \begin{aligned} & \ln \frac{\underline{\delta}_i+ |e_i\left(t\right)|/\rho_i\left(t\right)}{\bar{\delta}_i-|e_i\left(t\right)|/\rho_i\left(t\right)}, & \quad  \hspace{5pt} {\rm if} \hspace{5pt} e_i\left(t\right)\geq 0 \\
& - \ln \frac{\underline{\delta}_i+|e_i\left(t\right)|/\rho_i\left(t\right)}{\bar{\delta}_i-|e_i\left(t\right)|/\rho_i\left(t\right)} ,&  \hspace{5pt} {\rm if} \hspace{5pt} e_i\left(t\right)< 0\\ \end{aligned}  \right.\\
& \quad  {\rm with} \;\; \bar{\delta}_i > \underline{\delta}_i
\end{aligned}
\end{equation}

Thus, the transformed error can be expressed in more compact form as follows:

\begin{equation}
\label{eq:eq57b}
\begin{aligned}
\epsilon_i =  \frac{1}{2}{\rm sign}\left(e_i\left(t\right)/\rho_i\left(t\right)\right)\cdot & \ln \left(\frac{\underline{\delta}_i+ |e_i\left(t\right)|/\rho_i\left(t\right)}{\bar{\delta}_i-|e_i\left(t\right)|/\rho_i\left(t\right)}\right),\\
 & {\rm with} \;\;\quad \bar{\delta}_i > \underline{\delta}_i
\end{aligned}
\end{equation}
And to attenuate the effect of chattering, the following form of the transformed error will be considered
\begin{equation}
\label{eq:eq58b}
\begin{aligned}
\epsilon_i =  \frac{1}{2\sqrt{\pi}}erf\left(\frac{\xi e_i\left(t\right)}{\rho_i\left(t\right)}\right)\cdot & \ln \left(\frac{\underline{\delta}_i+ |e_i\left(t\right)|/\rho_i\left(t\right)}{\bar{\delta}_i-|e_i\left(t\right)|/\rho_i\left(t\right)}\right),\\
 & {\rm with} \;\; \quad \bar{\delta}_i > \underline{\delta}_i
\end{aligned}
\end{equation}
where $erf(\xi e/\rho)=\frac{2}{\sqrt{\pi}}\int\limits_{0}^{\xi e}e^{-a^2}da$. $\xi>0$ is a design parameter.
\begin{rem} The primary role of $\xi$ is to make the $erf(\xi e)$ as close as possible to ${\rm sign}\left(e\right)$. Ideally, $\xi$ is selected to be as big as possible. For instance, $|erf(\xi e)|\approxeq 1$ when $ |e|>\Delta=\frac{2}{\xi}$. Therefore, if $\xi=200$ then $|erf\left(e\right)|\approxeq 1$ when $|e|>0.01$.  However, while the derivative is smooth the more one selects a big $\xi$ the more there is a risk of chattering.
\end{rem}
For simplification, let $x\left(t\right)=x$, $e\left(t\right)=e$, $\epsilon\left(t\right)=\epsilon$ and $\rho\left(t\right)=\rho$. After algebraic manipulations, the derivative of transformed error when $|e|/\rho \geq \Delta/\xi$ can be approximated by:
\begin{equation}
\label{eq:eq17}
\dot{\epsilon}_i = \frac{1}{2\rho_i}\left(\frac{1}{\underline{\delta}_i+|e_i|/\rho_i} + \frac{1}{\bar{\delta}_i - |e_i|/\rho_i}\right)\left(\dot{e}_i - \frac{e_i\dot{\rho}_i}{\rho_i}\right)
\end{equation}
\begin{rem}
	As mentioned earlier, the selection of the high gain $\xi$ can make the absolute value of the error function converge to 1 for a small value of the ratio $|e\left(t\right)|/\rho=\Delta/\rho$. In our analysis, we will use (\ref{eq:eq17}) to show that the control will generate a UUB error dynamic that will converge to a ball around zero with a radius that can be made as small as desired depending on the selection of $\xi$. Thus, the error may not converge to zero.
\end{rem} Let
\begin{equation}
\label{eq:eq17c}
r_i = \frac{1}{2\rho_i}\left(\frac{1}{\underline{\delta}_i+|e_i|/\rho_i} + \frac{1}{\bar{\delta}_i - |e_i|/\rho_i}\right)
\end{equation}
From \eqref{eq:eq7} and \eqref{eq:eq17}, the global synchronization of the transformed error can be obtained as
\begin{equation}\label{eq:eq18}
\begin{aligned}
\dot{\epsilon} = & R\left(L+B\right)(A_mx+f\left(x\right)+B_mu+w-\underline{f}(x_{0},t))\\
 &- R\dot{\varUpsilon}\,\varUpsilon^{-1}\,e\left(t\right)
\end{aligned}
\end{equation}
where the control at the level of each node is of the form $u_i = -c \epsilon_i + \nu$; the value of $\nu$  represents the part of the control action necessary to tackle the uncertainties and takes into account the estimation errors in the adaptation rule (see \eqref{eq:eq27});
$\epsilon=[\epsilon_1,\ldots,\epsilon_N]^{\top}\in \mathbb{R}^N $, $\varUpsilon={\rm diag}\{\rho_i\left(t\right)\}$ and $\dot{\varUpsilon}={\rm diag}\{\dot{\rho_i}\left(t\right)\}$, $i=1, \ldots,N$;
$R$ is such that $R = diag[r_1\left(t\right),\ldots,r_N\left(t\right)]$
with  $R>0$ and $\dot{R}<0$;  $\,\dot{\varUpsilon}\,\varUpsilon^{-1}<0$ with $lim_{t \to \infty}\,\dot{\varUpsilon}\,\varUpsilon^{-1} = 0$.
Before proceeding further, the following definitions are needed (see \cite{das_distributed_2010}).
\begin{defn}
	The global error $e\left(t\right)\in \mathbb{R}^N $ is uniformly ultimately bounded (UUB)  if there exists a compact set $\Omega \subset \mathbb{R}^N $ so that $\forall e\left(t_{0}\right) \in \Omega$ there exists a bound $B$ and a time $t_{f}(B,e\left(t_{0}\right))$, both independent of $t_0 \geq 0$, such that $||e\left(t\right)|| \leq B$ so that $\forall t > t_0+t_{f}$.
\end{defn}
\begin{defn}
	The control node trajectory $x_{0}\left(t\right)$ given by \eqref{eq:eq1} is cooperative UUB with respect to solutions of node dynamics \eqref{eq:eq3} if there exists a compact set $\Omega \subset \mathbb{R}^N $ so that $\forall x_{i}\left(t_{0}\right) - x_{0}\left(t_{0}\right) \in \Omega$, there exist a bound $B$ and a time $t_{f}\left(B, x\left(t_{0}\right) - x_{0}\left(t_{0}\right) \right)$, both independent of $t_0\geq 0$, such that $||x\left(t_{0}\right) - x_{0}\left(t_{0}\right)|| \leq B$, $\forall i, \forall t > t_0+t_{f}$.
\end{defn}
\section{Adaptive Projection Approximation}\label{sec_4}
Using linear parametrization of nonlinear systems (for more details see A.8 in \cite{hovakimyan_l1_2010} )The agent $i$'s nonlinear dynamics in \eqref{eq:eq1} can be written as
\begin{equation}
\label{eq:eq19}
\begin{aligned}
& \dot{x}_i = A_{mi}x_{i} + B_{mi}u_i + \theta_i||x_{i}||_{\infty} + \sigma_i(x_{i},t)
\end{aligned}
\end{equation}
with $\theta_i \in \mathbb{R}$ is an unknown but bounded time varying parameter and $\sigma_i \in \mathbb{R}$ is the part that includes all unknown nonlinearities and external disturbances, $w$. $\Theta_i$ and $\Delta_i$ are known compact sets where $\theta_i \in \Theta_i$ and $\sigma_i \in \Delta_i$. In the remaining of the paper, the following assumptions will be considered.\\

\begin{assum} {\cite{hovakimyan_l1_2010}}         \label{assump2}
	\begin{enumerate}
		\item Leader's states are bounded by $||x_{0}|| \leq x_{0}$.
		\item Leader's nonlinear dynamic is unknown and bounded such as $||\underline{f}_0(x_{0},t)|| \leq F_M$.
		\item Uniform boundedness of the unknown parameters:   $||\theta\left(t\right)|| \leq \theta_M$ and $||\sigma\left(t\right)|| \leq \sigma_M$ for all $t>0$
		\item Uniform boundedness of the rate of variation of parameters:
		$\theta\left(t\right)$ and $\sigma\left(t\right)$ are continuously differentiable with uniformly bounded derivatives. $||\dot{\theta}\left(t\right)|| \leq d_{\theta} < \infty $ and $||\dot{\sigma}\left(t\right)|| \leq d_{\sigma} < \infty$  for all $ t \geq 0.$
	\end{enumerate}
	
	One should note that the values of the estimation bounds are not necessary known.
	\begin{assum}
		\label{assump1}
		Matrix $A_{mi}, B_{mi}$ are known and $B_{mi}^{-1}$ exists.
	\end{assum}
	Let $\hat{\theta}$ and $\hat{\sigma}$ be the approximation of $\theta$ and $\sigma$ respectively. Then,
	\begin{equation}
	\label{eq:eq20}
	\tilde{\theta}_i = \theta_i - \hat{\theta}_i
	\end{equation}
	\begin{equation}
	\label{eq:eq20_1}
	\tilde{\sigma}_i = \sigma_i - \hat{\sigma}_i
	\end{equation}
	\begin{rem}
		The communication graph is considered strongly connected. Thus, if $b_i \neq 0$ for at least one $i$, $i=1,\ldots, N$ then $\left(L+B\right)$ is an irreducible diagonally dominant M-matrix and hence nonsingular \cite{das_distributed_2010}.         The control signal of local agent $i$ can be given by
	\end{rem}
\end{assum}
\begin{lem} (see \cite{lewis_cooperative_2013} for more details.)
	\label{lemma2}
	Let $L$ be an irreducible matrix and $B \neq 0 $ such as $\left(L+B\right)$ is nonsingular, then we can define
	\begin{equation}
	\label{eq:eq24}
	q = [q_1,\ldots,q_N]^{\top}=\left(L+B\right)^{-1}\cdot\underline{1}
	\end{equation}
	\begin{equation}
	\label{eq:eq25}
	P = {\rm diag}\{p_i\}={\rm diag}\{1/q_i\}
	\end{equation}
	Then, $P > 0$ and the matrix $Q$ defined as
	\begin{equation}
	\label{eq:eq26}
	\begin{aligned}
	Q = & P\,\left(L+B\right)+\left(L+B\right)^{\top}\,P\\
	  = & P\left[S\left(L+B\right)+\left(L+B\right)^{\top}\,S\right]P
	\end{aligned}
	\end{equation}
	is also positive definite with  $S=P^{-1}$
\end{lem}
The gist of the idea is that $Q=S\left(L+B\right)+\left(L+B\right)^{\top}\,S$ is diagonally strictly dominant, and since it is a symmetric M-matrix, then it is positive definite. Based on this lemma, the following Preposition holds
\begin{prop}
	Let $R$ a positive definite diagonal matrix, and $L$, $B$, $P$ and $S$ as defined in Lemma 1, then the matrix $Q$ defined as
	\begin{equation}
	\label{eq:eq261}
	Q = P\,R\,\left(L+B\right)+\left(L+B\right)^{\top}\,R\,P
	\end{equation}
	is positive definite.
\end{prop}
\textbf{Proof:}\\
Since $\left(L+B\right)$ is is a nonsingular M-matrix and $R>0$ is diagonal, then $R\left(L+B\right)$ is a non-singular M-Matrix.
\begin{equation}
\left(L+B\right)\,q ={\bf \underline{1}}>0
\end{equation}
Let $S= {\rm diag}\{q_i\}$ then
\begin{equation}
R\left(L+B\right)S {\bf \underline{1}} = R\,\left(L+B\right)\,q=R{\bf \underline{1}}>0
\end{equation}
which means strict diagonal dominance of $R\,\left(L+B\right)\,S$.
\begin{equation}
\begin{aligned}
Q = & P\,R\,\left(L+B\right)+\left(L+B\right)^{\top}\,R\,P\\
  = & P\left[R\left(L+B\right)S+S\,\left(L+B\right)^{\top}\,R\right]P
\end{aligned}
\end{equation}
$R\left(L+B\right)S+S\,\left(L+B\right)^{\top}\,R$ is symmetric and strictly diagonally dominant. Therefore,  $Q$ is positive definite.\\

The control signal of local nodes is given by
\begin{equation}
\label{eq:eq27}
u_i = B_{mi}^{-1}\left( -c\epsilon_i - A_{mi}\,\tilde{x}_i - \hat{\theta}_i||x_{i}||_{\infty} - \hat{\sigma}_i \right)
\end{equation}
where the control gain $c>0$ and the overall control signal
\begin{equation}
\label{eq:eq28b}
u = B_{m}^{-1}\left( -c\epsilon - A_{m}\,\tilde{x} - \hat{\theta}\left  \Vert x \right \Vert_{\infty} - \hat{\sigma} \right)
\end{equation}
with $\left  \Vert x \right \Vert_{\infty} = [||x_1||_{\infty},\ldots,||x_N||_{\infty}]^{\top}$. Let the adaptive estimates  of $\hat{\theta}$ and $\hat{\sigma}$ updated according to
\begin{equation}
\label{eq:eq29}
\dot{\hat{\theta}}_{i} = \left(\Gamma_{i} x_{i}\epsilon_i^{T}p_{i}r_{i}(d_{i}+b_{i})\right)^{\top} - k\Gamma_{i}\hat{\theta}_{i}
\end{equation}
\begin{equation}
\label{eq:eq29_1}
\dot{\hat{\sigma}}_{i} = \left(\Gamma_{i} \epsilon_i^{T}p_{i}r_{i}(d_{i}+b_{i})\right)^{\top} - k\Gamma_{i}\hat{\sigma}_{i}
\end{equation}
with $\Gamma_i \in \mathbb{R}^{+}$ and $k>0$. $c$ and $k$ are scalar design parameters.
\begin{thm}
	Consider the strong connected digraph of the network in \eqref{eq:eq1} with adaptive estimates in \eqref{eq:eq29} and \eqref{eq:eq29_1} satisfying  Assumptions (\ref{assump2} and under the control law \label{eq:eq28}, then the distributed multi-agent system is UUB stable if the tuning gain $k$ and $c$ satisfy the following conditions
	\begin{equation}
	\label{eq:eq30}
	k=\frac{c\underline{\lambda}\left(Q\right)}{2}
	\end{equation}
	and
	\begin{equation}
	\label{eq:eq31}
	c\underline{\lambda}\left(Q\right)> \frac{1}{2}(x_m+1)\bar{\lambda}(P)\bar{\lambda}(A)
	\end{equation}
	with $P$ and $Q$ are defined in Lemma \ref{lemma2}.
\end{thm}

\textbf{Proof:}
using \eqref{eq:eq28b}, equation \eqref{eq:eq7} becomes
\begin{equation}
\label{eq:eq7bb}
\begin{aligned}
& \dot{e} = \left(L+B\right)(-c\epsilon+\tilde{\theta} \left  \Vert x \right \Vert_{\infty} + \tilde{\sigma}-\underline{f}(x_{0},t))
\end{aligned}
\end{equation}
Consider the following Lyapunov candidate function
\begin{equation}
\label{eq:eq34}
V = \frac{1}{2}\epsilon^{\top}P\epsilon + \frac{1}{2}\tilde{\theta}^{\top}\Gamma^{-1}\tilde{\theta} + \frac{1}{2}\tilde{\sigma}^{\top}\Gamma^{-1}\tilde{\sigma}
\end{equation}
with $P>0$ as defined in \eqref{eq:eq25}, $\gamma \in \mathbb{R}^+$ was mentioned in \eqref{eq:eq29} and $\Gamma = {\rm diag}\{\gamma_i\}$. The derivative of \eqref{eq:eq34} is
\begin{equation}
\label{eq:eq35}
\dot{V} = \epsilon^{\top}\,P\,\dot{\epsilon} + \tilde{\theta}^{\top}\,\Gamma^{-1}\,\dot{\tilde{\theta}} + \tilde{\sigma}^{\top}\,\Gamma^{-1}\,\dot{\tilde{\sigma}}
\end{equation}
Let $P_1=P\,R=R\,P$ and $Q=P_1\,\left(L+B\right)\,+\,\left(L+B\right)^{\top}\,P_1$. Using equations \eqref{eq:eq28b}-\eqref{eq:eq30} and \eqref{eq:eq7bb} to replace $\dot{\epsilon}$, $\dot{\tilde{\theta}}$ and $\dot{\tilde{\sigma}}$ respectively,  one can write
\begin{equation}
\label{eq:eq38}
\begin{split}
\dot{V} = & -\frac{1}{2}c\epsilon^{\top}Q\epsilon - \epsilon^{\top}P_1\left(L+B\right)\underline{f}(x_{0},t)) - k\tilde{\theta}^{\top}\tilde{\theta} - k\tilde{\sigma}^{\top}\tilde{\sigma}  \\
& + k\tilde{\sigma}^{\top}\sigma + k \tilde{\theta}^{\top}\theta + \epsilon^{\top}P_1A\tilde{\theta}\left  \Vert x \right \Vert_{\infty} + \epsilon^{\top}P_1A\tilde{\sigma}\\
&  - \epsilon^{\top}\,P_1\,\,\dot{\varUpsilon}\,\varUpsilon^{-1}\,e\left(t\right) +\tilde{\theta}^{\top}\,\Gamma^{-1}\,\dot{\theta} + \tilde{\sigma}^{\top}\,\Gamma^{-1}\,\dot{\sigma}
\end{split}
\end{equation}
On the other hand,
\begin{equation*}
e\left(t\right)= \varUpsilon\,S(\epsilon)
\end{equation*}
Therefore
\begin{equation}
\label{eq:eq38b}
\begin{split}
\dot{V} = & -\frac{1}{2}c\epsilon^{\top}Q\epsilon - \epsilon^{\top}P_1\left(L+B\right)\underline{f}(x_{0},t)) - k\tilde{\theta}^{\top}\tilde{\theta} - k\tilde{\sigma}^{\top}\tilde{\sigma} \\
&  + k\tilde{\sigma}^{\top}\sigma + k \tilde{\theta}^{\top}\theta + \epsilon^{\top}P_1A\tilde{\theta}\left  \Vert x \right \Vert_{\infty} + \epsilon^{\top}P_1A\tilde{\sigma} \\
&  - \epsilon^{\top}\,P\,\,\dot{\varUpsilon}\,\,S(\epsilon) +\tilde{\theta}^{\top}\,\Gamma^{-1}\,\dot{\theta} + \tilde{\sigma}^{\top}\,\Gamma^{-1}\,\dot{\sigma}
\end{split}
\end{equation}
one should note that $\varLambda\left(t\right)=-P\,\,\dot{\varUpsilon}$ is a positive definite diagonal matrix for $\forall\,t$ and $t\xrightarrow{\lim}\infty, \varLambda\left(t\right)=0$
\begin{equation}
\label{eq:eq38c}
\begin{split}
\dot{V} = & -\frac{1}{2}c\epsilon^{\top}Q\epsilon - \epsilon^{\top}P_1\left(L+B\right)\underline{f}(x_{0},t)) - k\tilde{\theta}^{\top}\tilde{\theta} - k\tilde{\sigma}^{\top}\tilde{\sigma} \\
& + k\tilde{\sigma}^{\top}\sigma  + k \tilde{\theta}^{\top}\theta + \epsilon^{\top}P_1A\tilde{\theta}\left  \Vert x \right \Vert_{\infty} + \epsilon^{\top}P_1A\tilde{\sigma} + \epsilon^{\top}\,\varLambda\,\,\bar{\delta} \\
& +\tilde{\theta}^{\top}\,\Gamma^{-1}\,\dot{\theta} + \tilde{\sigma}^{\top}\,\Gamma^{-1}\,\dot{\sigma}
\end{split}
\end{equation}
$\bar{\delta}={\rm max}\{\bar{\delta}_1, \,\ldots, \bar{\delta}_N\}$.
\begin{equation}
\label{eq:eq39}
\begin{split}
\dot{V} \leq & -\frac{1}{2}c\underline{\lambda}\left(Q\right)\left  \Vert \epsilon \right \Vert^2 + \bar{\lambda}\left(P_1\right)\bar{\lambda}\left(L+B\right)F_M\left  \Vert \epsilon \right \Vert - k\left  \Vert \tilde{\theta} \right \Vert^2   \\
& - k\left  \Vert \tilde{\sigma} \right \Vert^2 + \bar{\lambda}\left(P_1\right)\bar{\lambda}(A)x_M\left  \Vert \tilde{\theta} \right \Vert\left  \Vert \epsilon \right \Vert + \bar{\lambda}\left(P_1\right)\bar{\lambda}(A)\left  \Vert \epsilon \right \Vert\left  \Vert \tilde{\sigma} \right \Vert \\
& + k\left  \Vert \tilde{\sigma} \right \Vert\sigma_M  + k\left  \Vert \tilde{\theta} \right \Vert\theta_M  + \bar{\delta}\bar{\lambda}\left(\varLambda\right)\left  \Vert \epsilon \right \Vert + \bar{\lambda}\left(\Gamma^{-1}\right)\,d_{\theta}\, \left  \Vert \tilde{\theta} \right \Vert \\
& + \bar{\lambda}\left(\Gamma^{-1}\right) \,d_{\sigma}\,\left  \Vert \tilde{\sigma} \right \Vert
\end{split}
\end{equation}
Define 
\begin{equation*}
\begin{aligned}
z = & \begin{bmatrix} \left  \Vert \epsilon \right \Vert & \left  \Vert \tilde{\theta} \right \Vert & \left  \Vert \tilde{\sigma} \right \Vert \end{bmatrix}^{\top}
\end{aligned}
\end{equation*}
{\small
\begin{equation*}
\begin{aligned}
H = & \begin{bmatrix}
\frac{1}{2}c\underline{\lambda}\left(Q\right)  & -\frac{1}{2}\bar{\lambda}\left(P_1\right)\bar{\lambda}(A)x_M & -\frac{1}{2}\bar{\lambda}\left(P_1\right)\bar{\lambda}(A)\\
-\frac{1}{2}\bar{\lambda}\left(P_1\right)\bar{\lambda}(A)x_M & k & 0\\
-\frac{1}{2}\bar{\lambda}\left(P_1\right)\bar{\lambda}(A) & 0 & k
\end{bmatrix}\\
h = & \begin{bmatrix}
\bar{\lambda}\left(P_1\right)\left( \bar{\lambda}\left(L+B\right)F_M + \bar{\delta}\bar{\lambda}\left(\varLambda\right) \right) \\
 k\theta_M+\bar{\lambda}\left(\Gamma^{-1}\right)\,d_{\theta}\\
 k\sigma_M+\bar{\lambda}\left(\Gamma^{-1}\right) \,d_{\sigma}
\end{bmatrix}
\end{aligned}
\end{equation*} }
then
	\begin{equation}
	\label{eq:eq40}
	\begin{split}
	\dot{V} \leq
	& -z^{\top}Hz +h^{\top}z
	\end{split}
	\end{equation}
 then equation \eqref{eq:eq40} can be written as
\begin{equation}
\label{eq:eq41}
\begin{split}
\dot{V} \leq &-z^{\top}Hz + h^{\top}z
\end{split}
\end{equation}
and $\dot{V} \leq 0$ if and only if $H$ is positive definite and
\begin{equation}
\label{eq:eq42}
\begin{split}
\left  \Vert z  \right \Vert > & \frac{\left  \Vert h \right \Vert}{\underline{\lambda}\left(H\right)}
\end{split}
\end{equation}
The Lyapunov candidate in \eqref{eq:eq34} can be described by
\begin{equation}
\label{eq:eq43}
\begin{split}
& \frac{1}{2}\underline{\lambda}(P)\left  \Vert \epsilon \right \Vert^2 + \frac{\underline{\lambda}\left(\Gamma^{-1}\right)}{2}\left  \Vert \tilde{\theta} \right \Vert^2 + \frac{\underline{\lambda}\left(\Gamma^{-1}\right)}{2}\left  \Vert \tilde{\sigma} \right \Vert^2 \\
& \hspace{40pt}\leq V \leq \\
& \frac{1}{2}\bar{\lambda}(P)\left  \Vert \epsilon \right \Vert^2  + \frac{\bar{\lambda}\left(\Gamma^{-1}\right)}{2}\left  \Vert \tilde{\theta} \right \Vert^2 + \frac{\bar{\lambda}\left(\Gamma^{-1}\right)}{2}\left  \Vert \tilde{\sigma} \right \Vert^2
\end{split}
\end{equation}
or
\begin{equation}
\label{eq:eq44}
\begin{split}
& \frac{1}{2}z^{\top}
\begin{bmatrix}
\underline{\lambda}(P)  & 0 & 0\\
0 & \underline{\lambda}\left(\Gamma^{-1}\right) & 0\\
0 & 0 & \underline{\lambda}\left(\Gamma^{-1}\right)
\end{bmatrix}
z \leq V \leq \\
&
\hspace{60pt}\frac{1}{2}z^{\top}
\begin{bmatrix}
\bar{\lambda}(P)  & 0 & 0\\
0 & \bar{\lambda}\left(\Gamma^{-1}\right) & 0\\
0 & 0 & \bar{\lambda}\left(\Gamma^{-1}\right)
\end{bmatrix}
z
\end{split}
\end{equation}
Let
\begin{equation*}
\Pi_{min}=\begin{bmatrix}
\underline{\lambda}(P)  & 0 & 0\\
0 & \underline{\lambda}\left(\Gamma^{-1}\right) & 0\\
0 & 0 & \underline{\lambda}\left(\Gamma^{-1}\right)
\end{bmatrix}
\end{equation*}
\begin{equation*}
\Pi_{max}= \begin{bmatrix}
\bar{\lambda}(P)  & 0 & 0\\
0 & \bar{\lambda}\left(\Gamma^{-1}\right) & 0\\
0 & 0 & \bar{\lambda}\left(\Gamma^{-1}\right)
\end{bmatrix}
\end{equation*}
\eqref{eq:eq44} is equivalent to
\begin{equation}
\label{eq:eq45}
\begin{split}
\frac{1}{2}\underline{\lambda}(\Pi_{min})\left  \Vert z  \right \Vert^2 \leq V \leq \frac{1}{2}\bar{\lambda}(\Pi_{max})\left  \Vert z  \right \Vert^2
\end{split}
\end{equation}
then
\begin{equation}
\label{eq:eq46}
\begin{split}
V > \frac{1}{2}\underline{\lambda}(\Pi_{min})\frac{\left  \Vert h \right \Vert^2}{\underline{\lambda}\left(H\right)^2}
\end{split}
\end{equation}
Defining $c=\frac{2}{\underline{\lambda}\left(Q\right)}$, $\gamma_1 = \frac{1}{2}\bar{\lambda}\left(P_1\right)\bar{\lambda}(A)x_M$, $\gamma_2 = \frac{1}{2}\bar{\lambda}\left(P_1\right)\bar{\lambda}(A)$ and
substitute \eqref{eq:eq40}, we have
\begin{equation*}
H = \begin{bmatrix}
c  & -\gamma_1   & -\gamma_2 \\
-\gamma_1 & k & 0 \\
-\gamma_2 & 0 & k \end{bmatrix}
\end{equation*}
where $H$ is positive definite matrix. If we select $k=\frac{1}{2}c\underline{\lambda}\left(Q\right)$ and $c\underline{\lambda}\left(Q\right)> \frac{1}{2}(x_M+1)\bar{\lambda}(P)\bar{\lambda}(A)$, then, we will have
\begin{equation*}
\underline{\lambda}\left(H\right) = \frac{c\underline{\lambda}\left(Q\right) - \frac{1}{2}(x_m+1)\bar{\lambda}(P)\bar{\lambda}(A)}{2}
\end{equation*}
And from \ref{eq:eq42}, taking $||\cdot||_1$ of $h$, define
\begin{equation*}
s_{1} = \bar{\lambda}\left(P_1\right) \left( \bar{\lambda}\left(L+B\right)F_M +  \bar{\delta}\bar{\lambda}\left(\varLambda\right) \right)
\end{equation*}
As such
	\begin{equation}
	\label{eq:eq48}
	\begin{split}
	\left  \Vert z  \right \Vert > \frac{s_{1}\bar{\lambda}(\dot{\rho})  +  k\theta_M+\bar{\lambda}\left(\Gamma^{-1}\right)\,d_{\theta} + k\sigma_M+\bar{\lambda}\left(\Gamma^{-1}\right) \,d_{\sigma}}{\underline{\lambda}\left(H\right)}
	\end{split}
	\end{equation}
which implies
	\begin{equation}
	\label{eq:eq49}
	\begin{split}
	\left  \Vert \epsilon \right \Vert > \frac{s_{1} + k\theta_M+\bar{\lambda}\left(\Gamma^{-1}\right)\,d_{\theta} + k\sigma_M+\bar{\lambda}\left(\Gamma^{-1}\right) \,d_{\sigma}}{\underline{\lambda}\left(H\right)}
	\end{split}
	\end{equation}
	\begin{equation}
	\label{eq:eq50}
	\begin{split}
	\left  \Vert \tilde{\theta} \right \Vert > \frac{s_{1} + k\theta_M+\bar{\lambda}\left(\Gamma^{-1}\right)\,d_{\theta} + k\sigma_M+\bar{\lambda}\left(\Gamma^{-1}\right) \,d_{\sigma}}{\underline{\lambda}\left(H\right)}
	\end{split}
	\end{equation}
	\begin{equation}
	\label{eq:eq50_1}
	\begin{split}
	\left  \Vert \tilde{\sigma} \right \Vert > \frac{s_{1} + k\theta_M+\bar{\lambda}\left(\Gamma^{-1}\right)\,d_{\theta} + k\sigma_M+\bar{\lambda}\left(\Gamma^{-1}\right) \,d_{\sigma}}{\underline{\lambda}\left(H\right)}
	\end{split}
	\end{equation}
Also from \eqref{eq:eq45}, we have
\begin{equation}
\label{eq:eq51}
\begin{split}
\left  \Vert z  \right \Vert \leq \sqrt{\frac{2V}{\underline{\lambda}(\underline{S})}}, \hspace{10pt}\left  \Vert z  \right \Vert \geq \sqrt{\frac{2V}{\bar{\lambda}(\bar{S})}}
\end{split}
\end{equation}
Then, equation \eqref{eq:eq41} can be written as
\begin{equation}
\label{eq:eq52}
\dot{V} \leq  -H_{\alpha}V+ h_{\alpha}\sqrt{V}
\end{equation}
with $H_{\alpha} = \frac{2\underline{\lambda}\left(H\right)}{\bar{\lambda}(\bar{\Pi})}$ and $h_{\alpha} = \frac{\sqrt{2}\left  \Vert h \right \Vert}{\sqrt{\underline{\lambda}(\underline{\Pi})}}$ which equivalent to
\begin{equation}
\label{eq:eq53}
\frac{2 d\sqrt{V}}{dt} \leq  -H_{\alpha}\sqrt{V} + h_{\alpha}
\end{equation}
\begin{equation}
\label{eq:eq54}
\sqrt{V} \leq exp^{-H_{\alpha}\,t/2}\left(\sqrt{V(0)} - \frac{h_{\alpha}}{H_{\alpha}}\right) + \frac{h_{\alpha}}{H_{\alpha}}
\end{equation}
That can be written as
\begin{equation}
\label{eq:eq60}
\sqrt{V} \leq \sqrt{V(0)} \leq  \sqrt{V(0)} + \frac{h_{\alpha}}{H_{\alpha}}
\end{equation}

Finally, the algorithm of nonlinear single node dynamics such as equation \eqref{eq:eq1} can be summarized briefly as
\begin{enumerate}[1.]
	\item Define the system known parameters $A_{mi}, B_{mi}$.
	\item Define the control design parameters such as $\Gamma_i$, $p_i$, $d_i$, $b_i$, $k$ and $c$.
	\item Evaluate local error synchronization from equation  \eqref{eq:eq4}.
	\item Evaluate the prescribed performance function from equation \eqref{eq:eq9}.
	\item Evaluate $r_i$ from equation \eqref{eq:eq17c}.
	\item Evaluate transformed error from equation  \eqref{eq:eq58b}.
	\item Evaluate control signal from equation \eqref{eq:eq27}.
	\item Evaluate adaptive estimates from equations \eqref{eq:eq29} and \eqref{eq:eq29_1}.
	\item Go to step 3).
\end{enumerate}

\begin{rem}
	If we have $x_{i} \in \mathbb{R}^{n}, n > 1$, then $u_i\in \mathbb{R}^{n}$, $f_i\left(x_{i}\right)\in \mathbb{R}^{n}$, $\theta_i \in \mathbb{R}^{n}$, $\sigma_i \in \mathbb{R}^{n}$, $r_{i} = {\rm diag}\{r_{i1},\ldots,r_{in}\} \in \mathbb{R}^{n \times n}$, then the problem can be extended easily and the estimated weight will be written as
	\begin{equation}
	\label{eq:eq55}
	\dot{\hat{\theta}}_{i} = \left(\Gamma_{i} ||x_{i}||_{\infty}\epsilon_i^{T}(p_{i} \mathbb{I}_n)r_{i}((d_i+b_i) \mathbb{I}_n)\right)^{\top} - k\Gamma_{i}\hat{\theta}_{i}
	\end{equation}
	\begin{equation}
	\label{eq:eq55_1}
	\dot{\hat{\sigma}}_{i} = \left(\Gamma_{i} \epsilon_i^{T}(p_{i} \mathbb{I}_n)r_{i}((d_i+b_i) \mathbb{I}_n)\right)^{\top} - k\Gamma_{i}\hat{\sigma}_{i}
	\end{equation}
\end{rem}

\section{Examples and Simulations}\label{sec_5}

\begin{figure}[h!]
	\centering
	\includegraphics[scale=0.5]{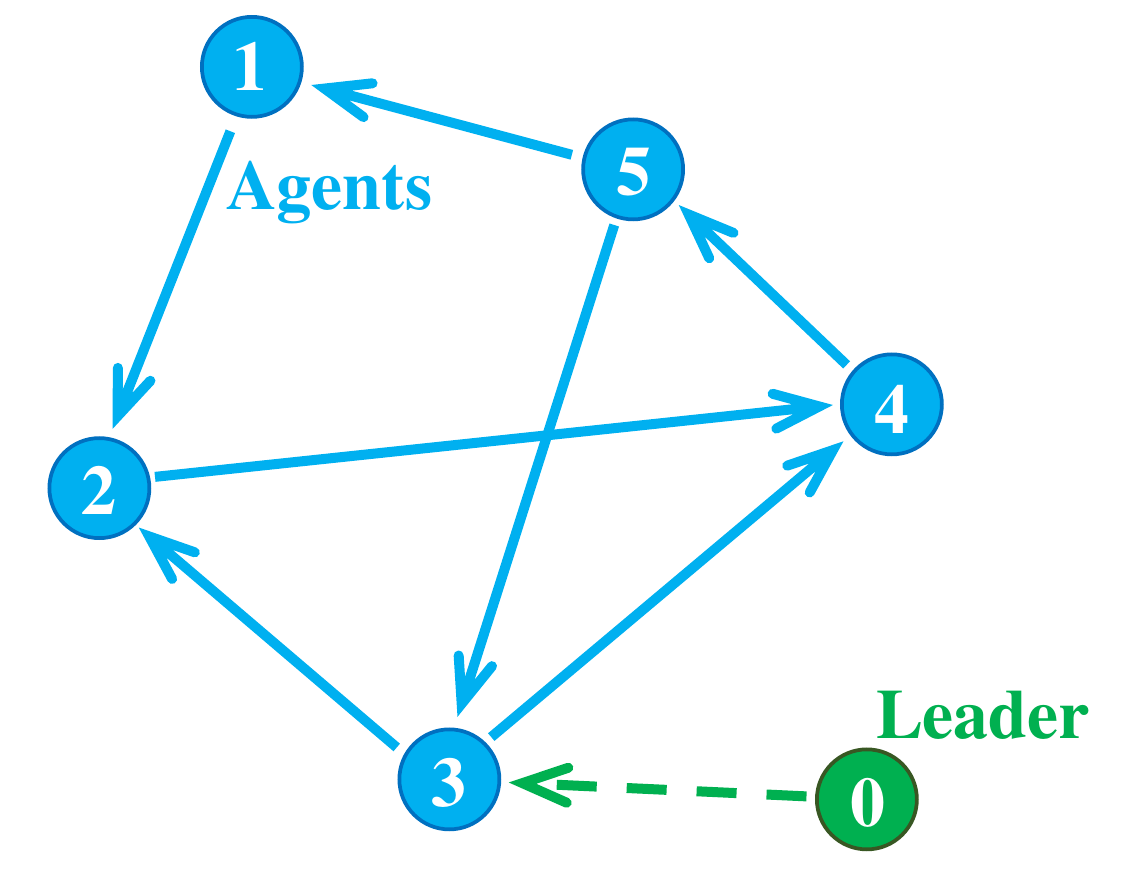}
	\caption{ Strongly connected graph of one leader and five agents.}
	\label{fig:fig2}
\end{figure}

\begin{figure*}[ht]
	\centering
	\includegraphics[scale=0.5]{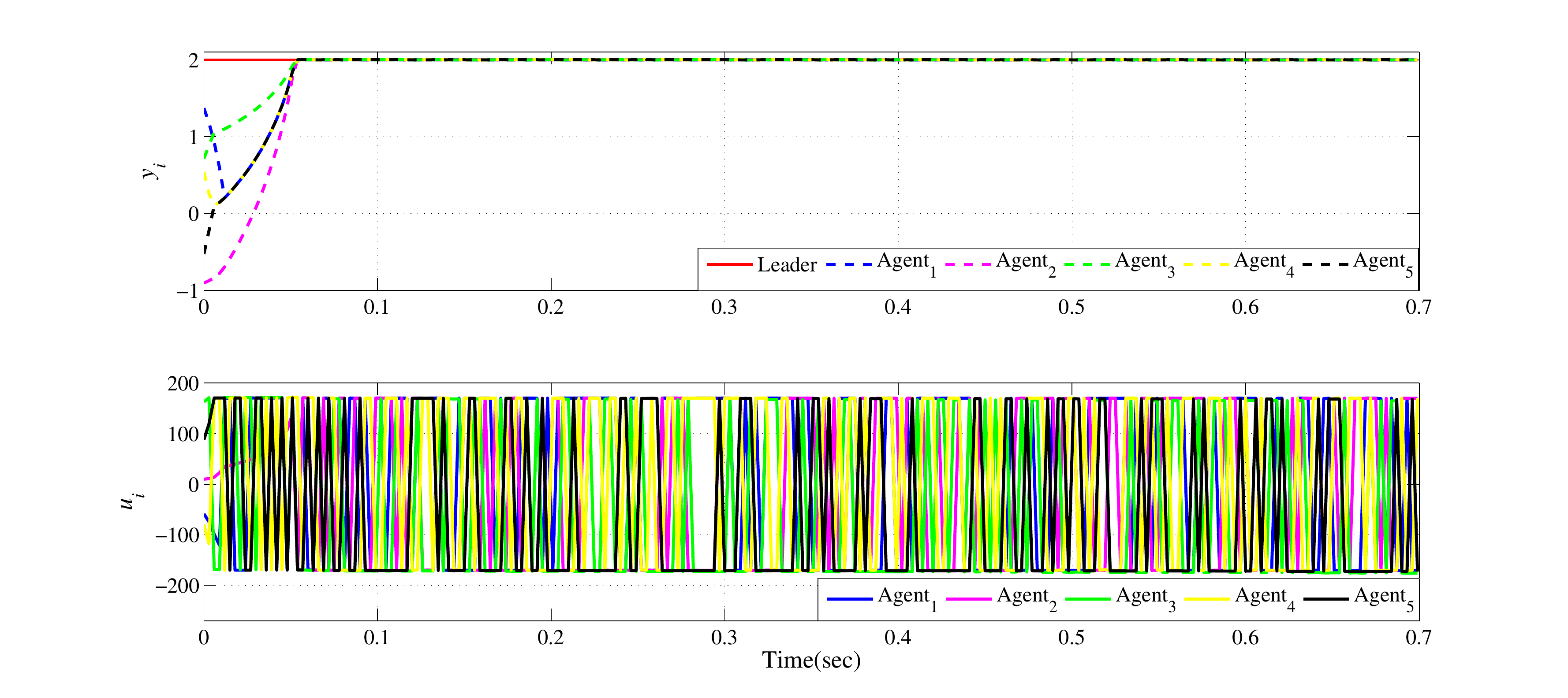}
	\caption{ Output performance and control signal using \eqref{eq:eq16b} for the transformed error.}
	\label{fig:fig3}
\end{figure*}

\begin{figure*}[ht]
	\centering
	\includegraphics[scale=0.5]{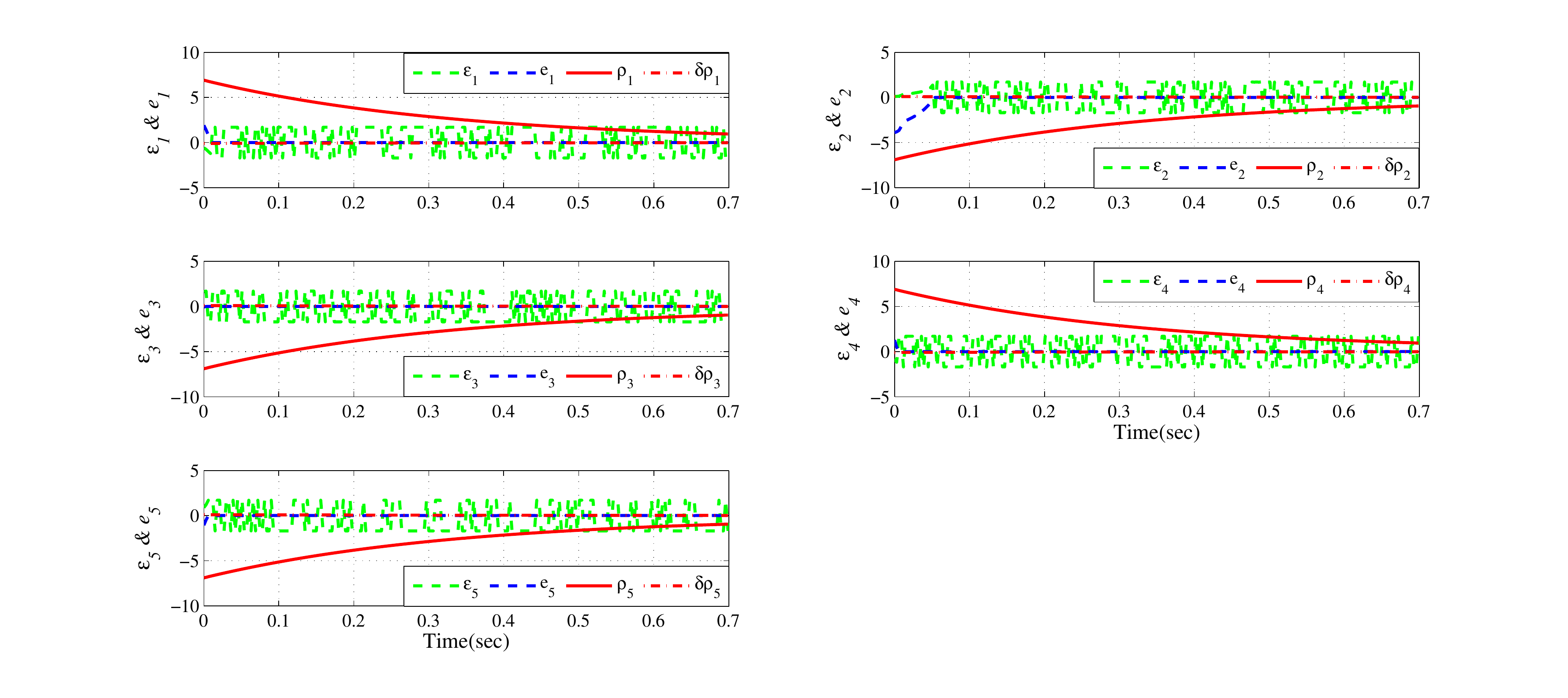}
	\caption{ Error and Transformed Error using \eqref{eq:eq16b} for the transformed error.}
	\label{fig:fig4}
\end{figure*}

{\bf Example 1:} Consider the digraph composed of five nodes strongly connected and having a single leader connected to node 3. The pining gains between connected nodes are assumed equal to 1 as in Fig. \ref{fig:fig2}.
The nonlinear dynamics of the different agents are as follows
\begin{equation*}
\begin{aligned}
& \dot{x}_1 = x_1^3 + u_1 + a_1 cos\left(t\right)\\
& \dot{x}_2 = x_2^2 + u_2 + a_2 cos\left(t\right)\\
& \dot{x}_3 = x_3^4 + u_3 + a_3 cos\left(t\right)\\
& \dot{x}_4 = x_4 + u_4 + a_4 cos\left(t\right)\\
& \dot{x}_5 = x_5^5 + u_5 + a_5 cos\left(t\right)
\end{aligned}
\end{equation*}
$a_i$, $i= 1,\ldots, 5$ are bounded randomly generated constant amplitudes. The leader dynamics was selected $\dot{x}_0 = f_0(x_{0},t) = 0$ with desired consensus value equal to $2$. Nonlinearities and disturbances are assumed to be unknown in all nodes. The control parameters of the system are
$\rho_{\infty} = 0.05\times{\bf \underline{1}_{1\times 5}}$, $\rho_{0} = 7\times{\bf \underline{1}_{1\times 5}}$, $l = 7\times{\bf \underline{1}_{1\times 5}}$, $\Gamma = 150\mathbb{I}_{5\times 5}$, $\bar{\delta} = 7$, $\underline{\delta} = 1$, $c = 100$, $k = 0.8$, $\alpha = 20$ and $x_{0} = 2$,
$x(0) = [0.8230, -0.9001, -2.5351,-1.4567, -0.7553]^{\top}$.\\
Figures \ref{fig:fig3} and \ref{fig:fig4} show the output performance, control signal and transformed error respectively for the proposed control algorithm using \eqref{eq:eq16b} for the transformed error. Fig. \ref{fig:fig3} shows the severe chattering in the control effort. Although the oscillation of synchronization error values satisfy the prescribed performance conditions, the switching in error signs caused switching in transformed errors as shown in Fig. \ref{fig:fig4} which consequently causes chattering in the control signal as clearly revealed in Fig. \ref{fig:fig3}.

\begin{figure*}[ht]
	\centering
	\includegraphics[scale=0.5]{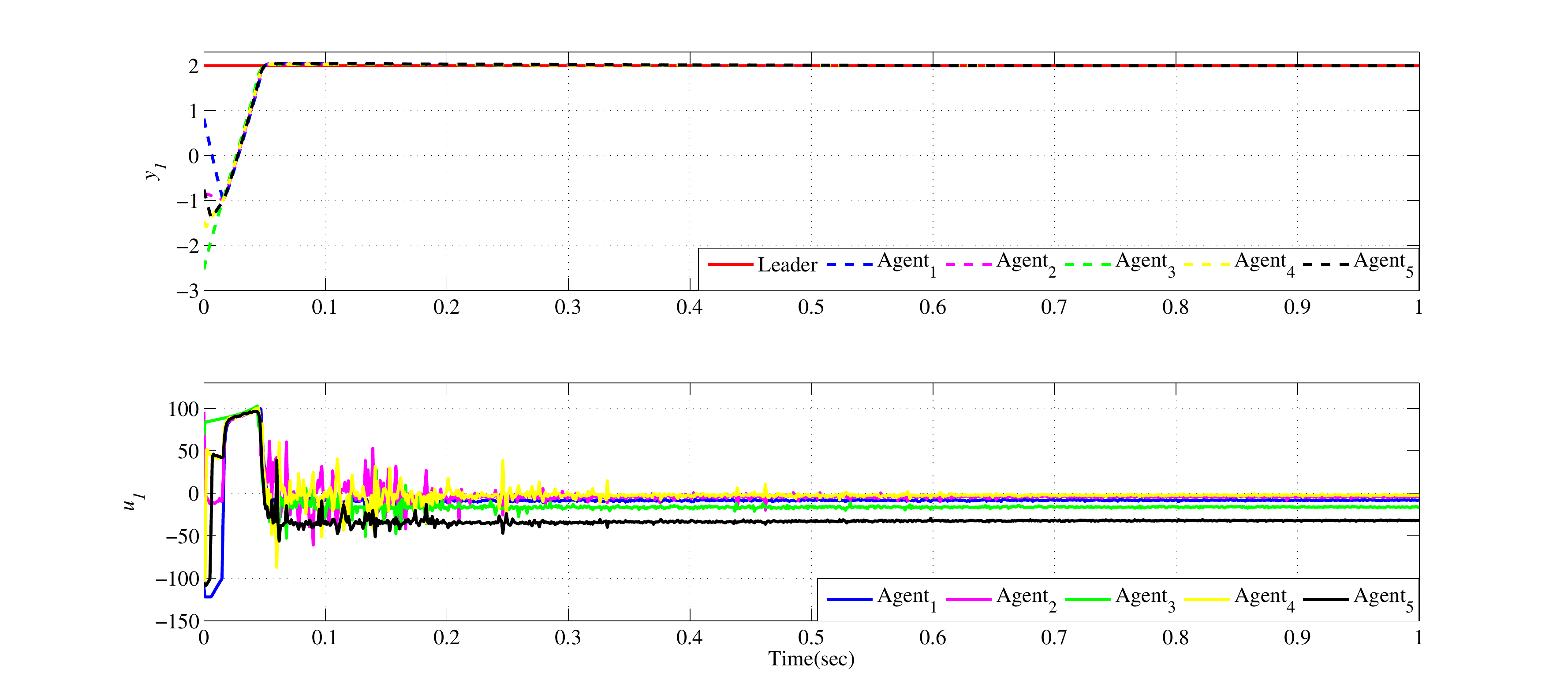}
	\caption{Output performance and control signal with the new prescribed performance function as in \eqref{eq:eq58b}.}
	\label{fig:fig5}
\end{figure*}

\begin{figure*}[ht]
	\centering
	\includegraphics[scale=0.5]{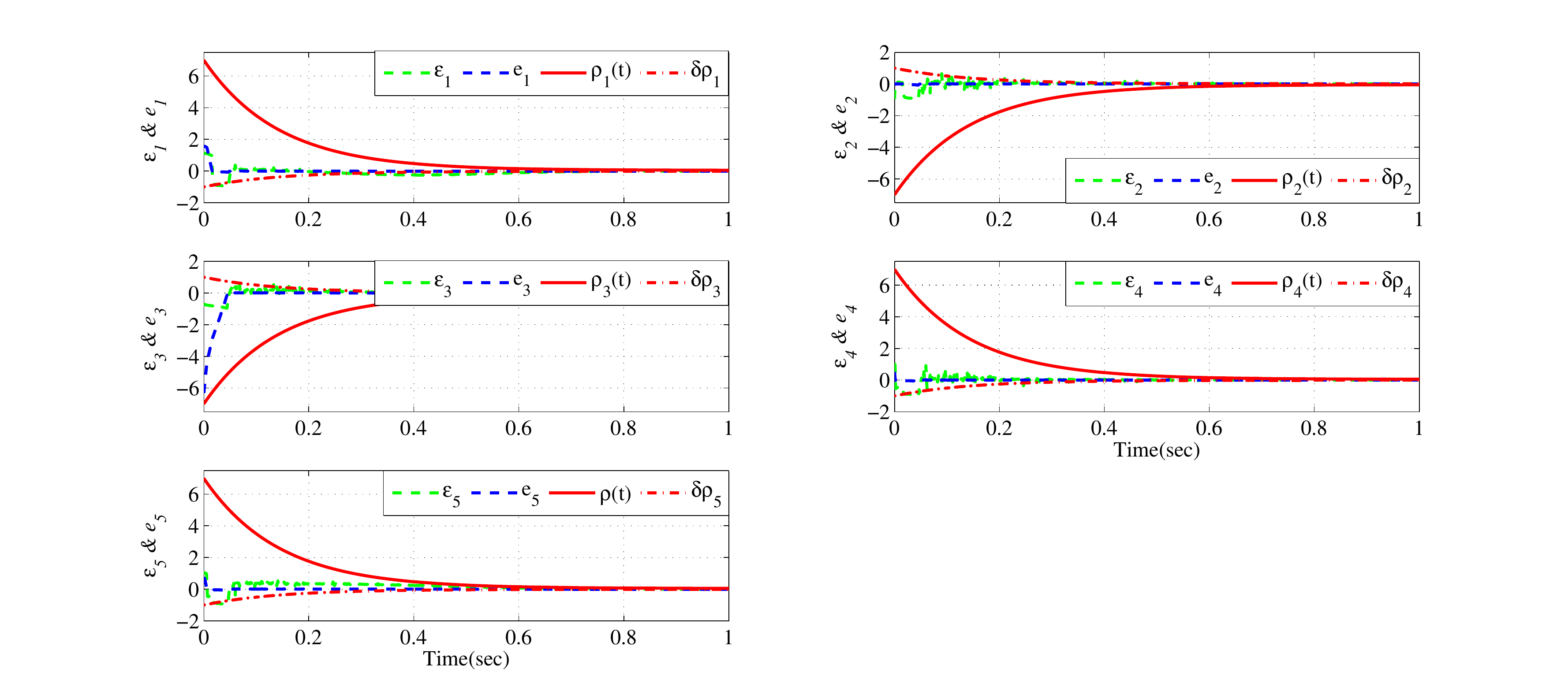}
	\caption{Error and transformed error with the new prescribed performance function as in \eqref{eq:eq58b}.}
	\label{fig:fig6}
\end{figure*}

The proposed control with the new prescribed performance function as in \eqref{eq:eq58b} Fig. \ref{fig:fig5} and \ref{fig:fig6} show the output performance, control signal and transformed error respectively. A significant improvement in the control effort and transformed error can be clearly observed.\\

%

\begin{figure*}[ht]
	\centering
	\includegraphics[scale=0.5]{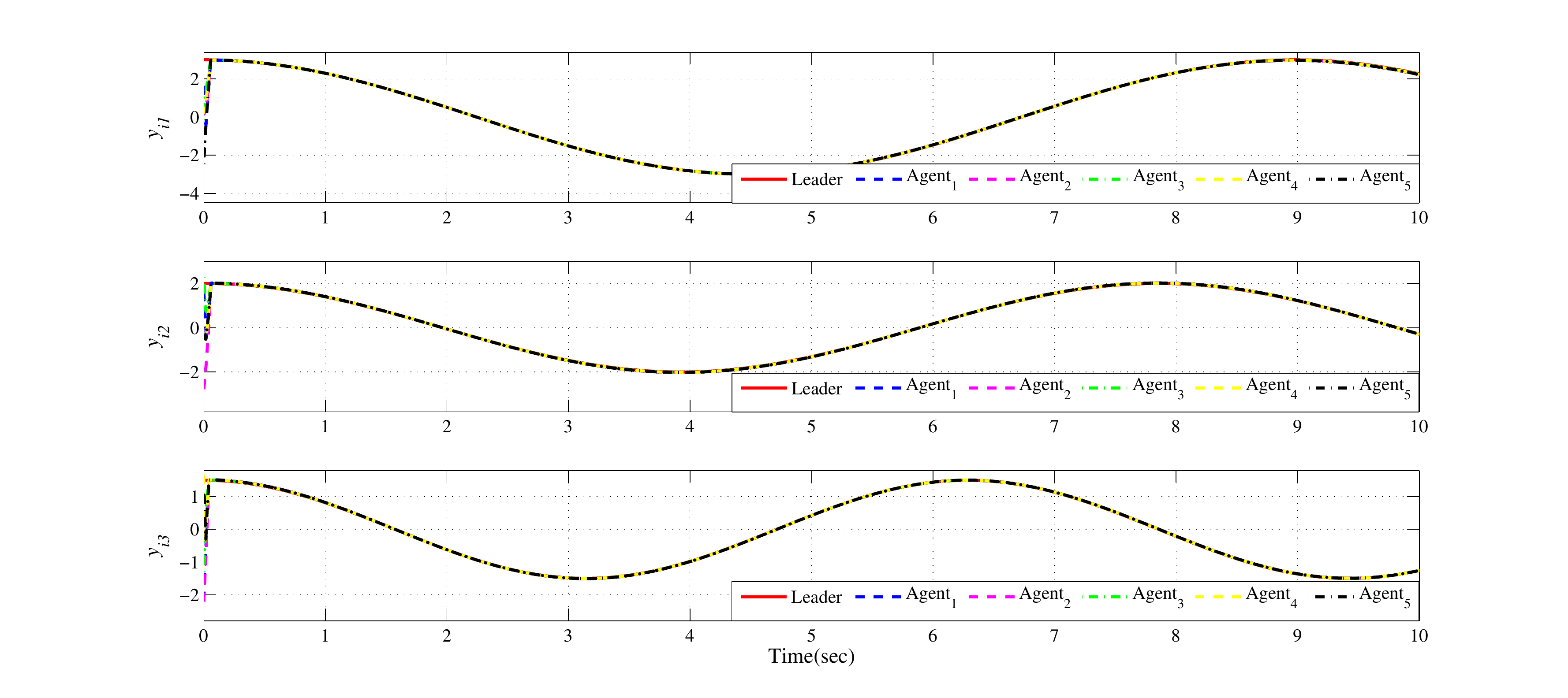}
	\caption{Output performance in the MIMO case using \eqref{eq:eq58b}.}
	\label{fig:fig9}
\end{figure*}

\begin{figure*}[ht]
	\centering
	\includegraphics[scale=0.5]{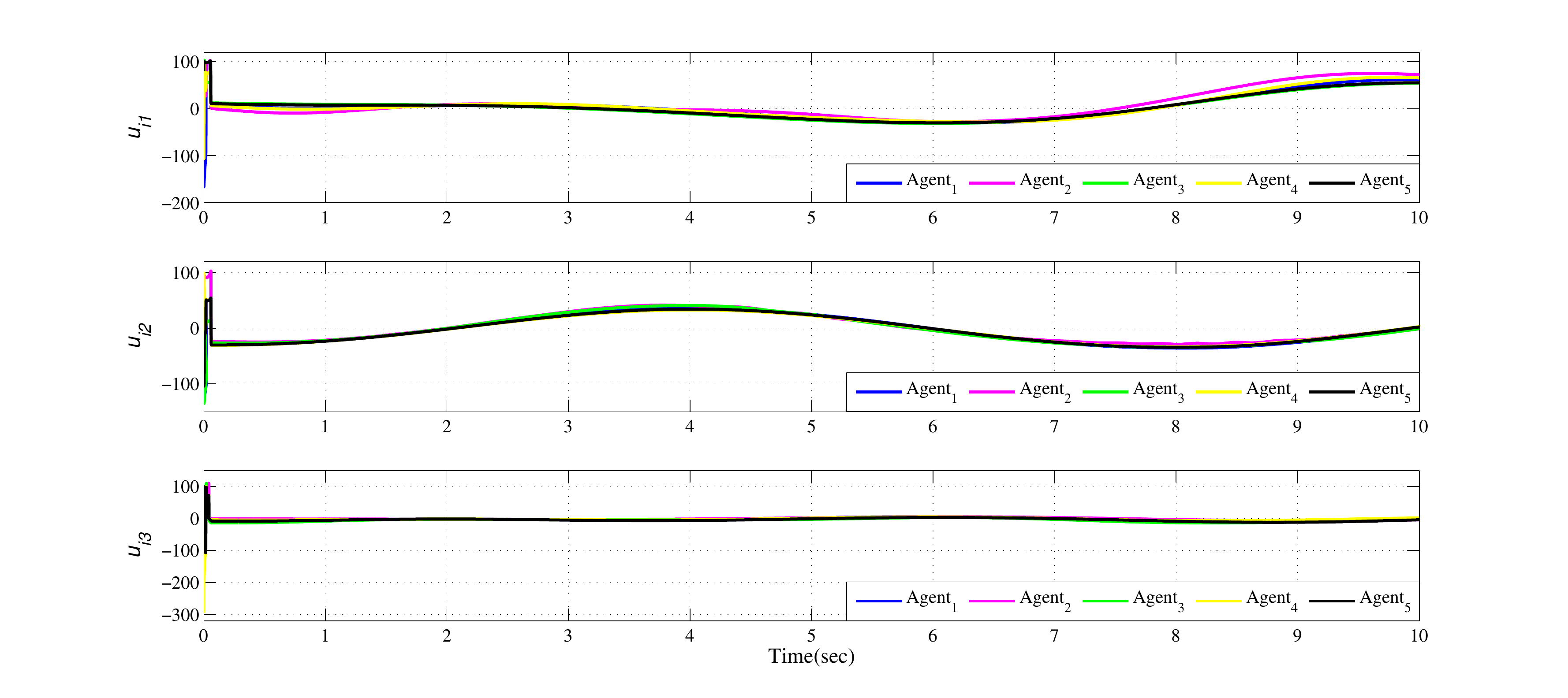}
	\caption{Control Signal in the MIMO case using \eqref{eq:eq58b}.}
	\label{fig:fig10}
\end{figure*}

\begin{figure*}[ht]
	\centering
	\includegraphics[scale=0.5]{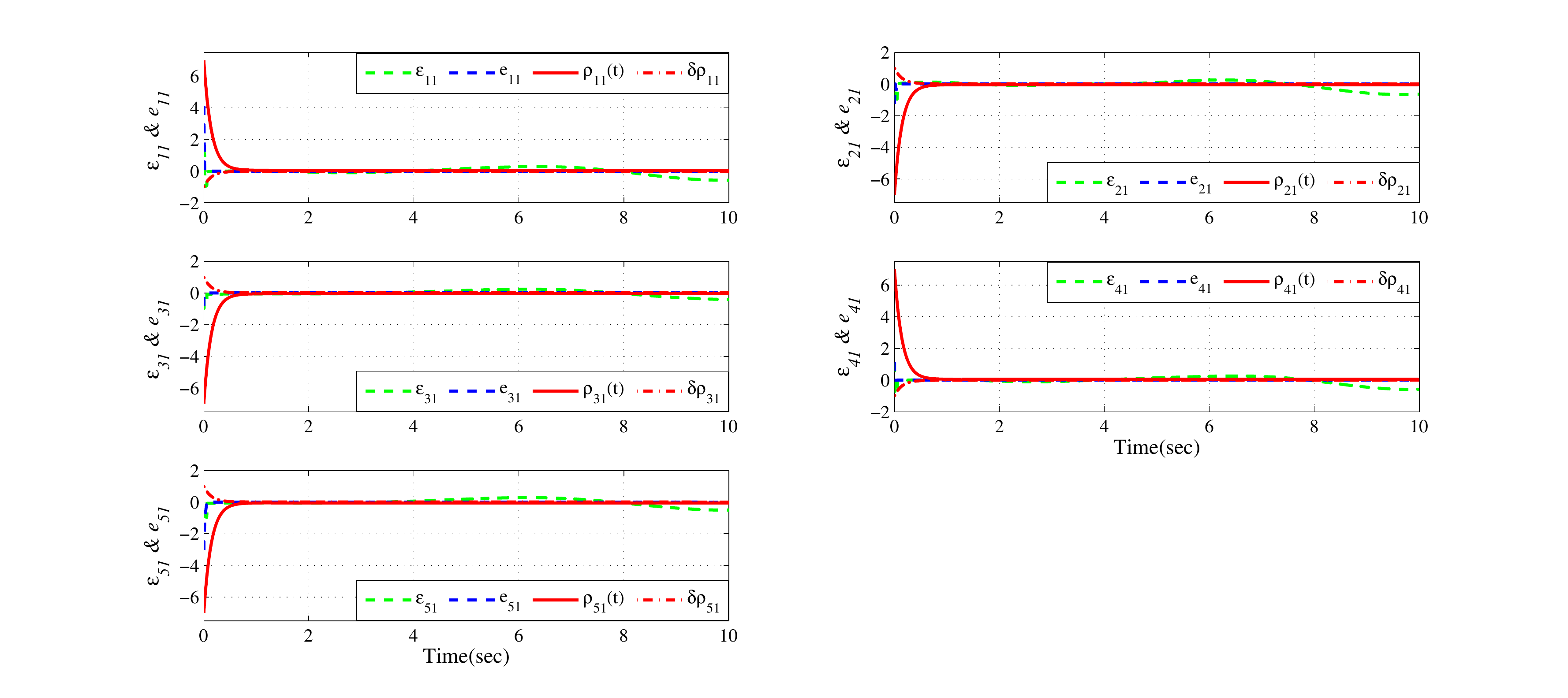}
	\caption{ Error and Transformed Error in the MIMO case using \eqref{eq:eq58b} for $x_{1,j}$ where $j = 1, \ldots, 5$.}
	\label{fig:fig11}
\end{figure*}

\begin{figure*}[ht]
	\centering
	\includegraphics[scale=0.5]{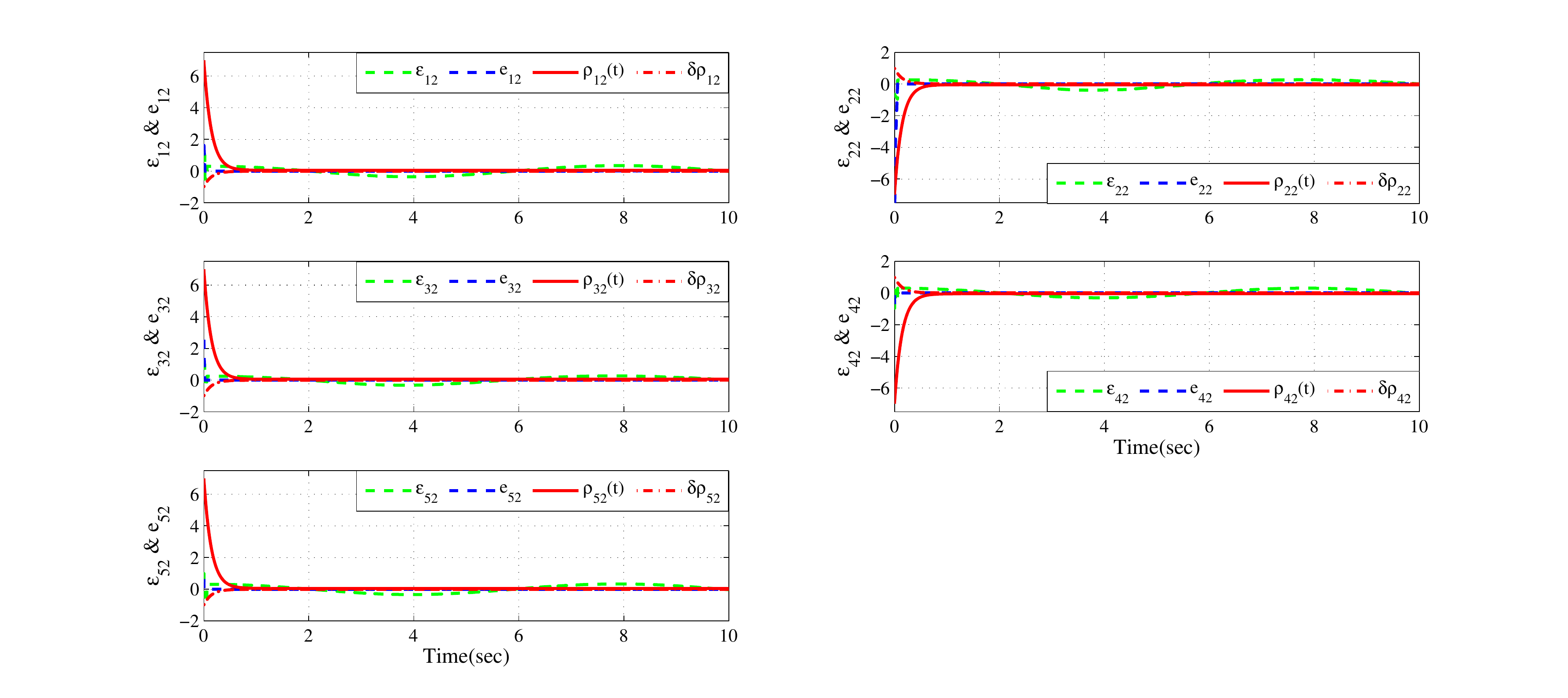}
	\caption{ Error and Transformed Error in the MIMO case using \eqref{eq:eq58b} for $x_{2,j}$ where $j = 1, \ldots, 5$.}
	\label{fig:fig12}
\end{figure*}

\begin{figure*}[ht]
	\centering
	\includegraphics[scale=0.5]{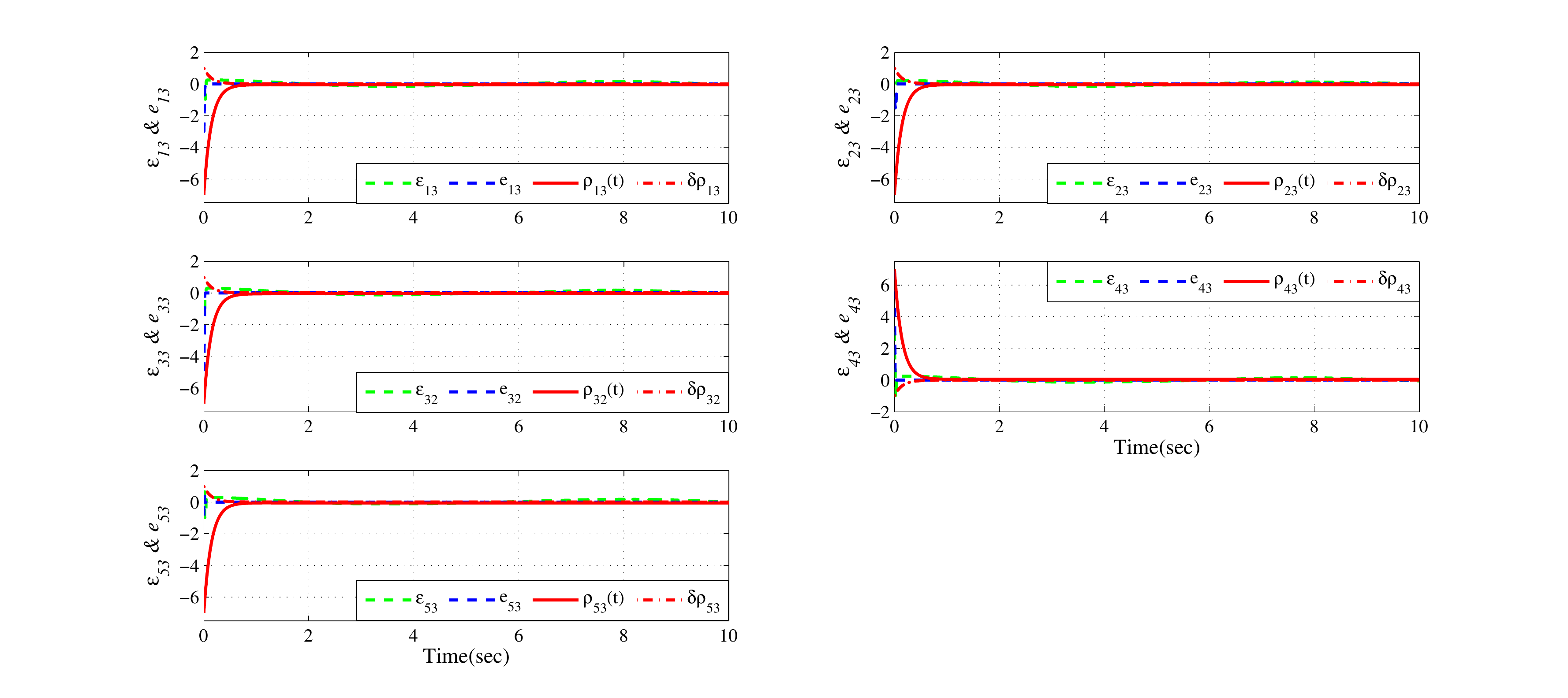}
	\caption{ Error and Transformed Error in the MIMO case using \eqref{eq:eq58b} for $x_{3,j}$ where $j = 1, \ldots, 5$}
	\label{fig:fig13}
\end{figure*}

\begin{figure*}[ht]
	\centering
	\includegraphics[scale=0.5]{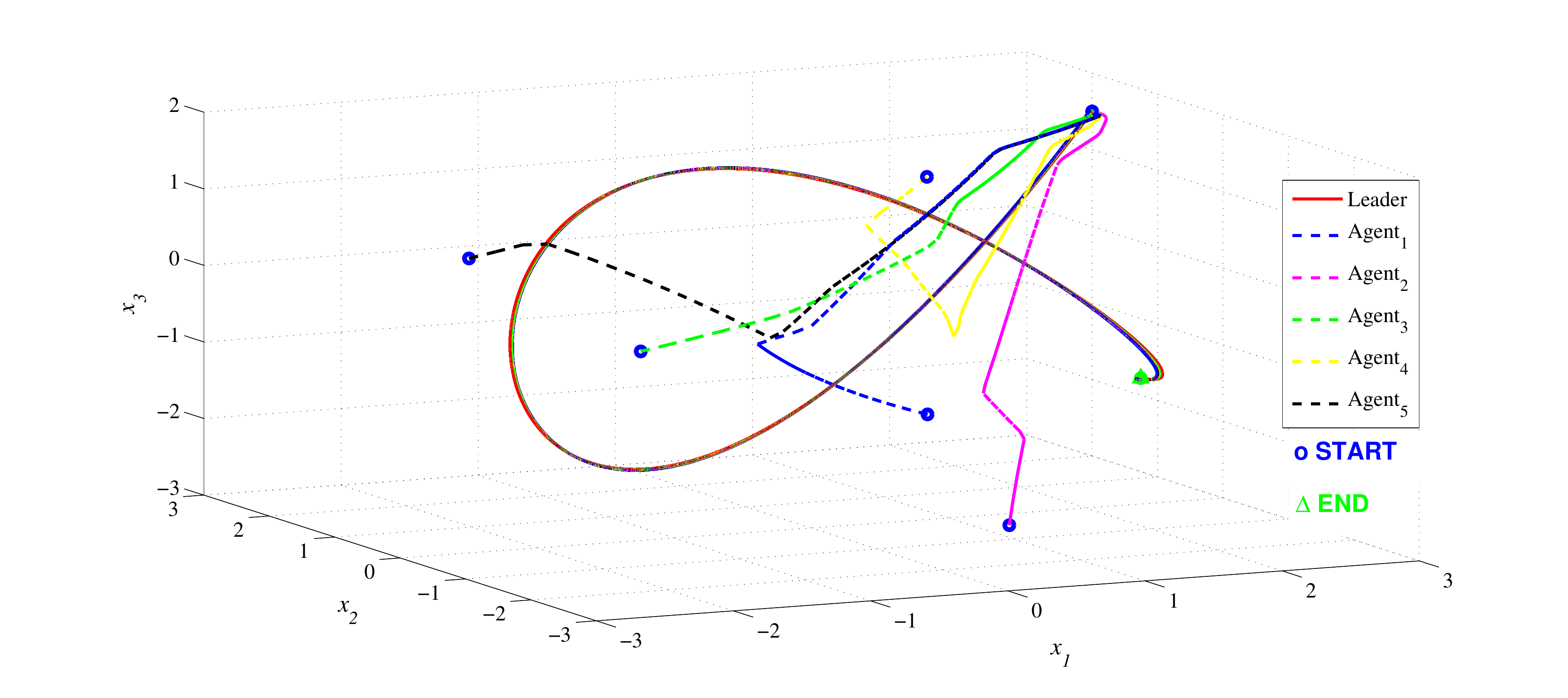}
	\caption{Phase Plan for motion synchronization in the MIMO system case under different initial conditions.}
	\label{fig:fig14}
\end{figure*}

{\bf Example 2: (MIMO case)} Consider the same problem as in Fig. \ref{fig:fig2} with 3 inputs and 3 outputs nonlinear systems. The nonlinear dynamics of the graph are now
\begin{equation*}
\begin{aligned}
& \dot{x}_{j} = Ax_{j} + Bu_{j} + \theta_j x_{j} + f_j(x_{j}) + D_j\left(t\right) \\
& y_{j} = Cx_{j}\\
\end{aligned}
\end{equation*}
where $x_{j} \in \mathbb{R}^{3 \times 1}$ is the state vector,  $u_{j} \in \mathbb{R}^{3 \times 1}$ is the input vector,  $y_{j} \in \mathbb{R}^{3 \times 1}$ is the output vector. $A, B and C \in \mathbb{R}^{3 \times 3}$ are known constant matrices.
\begin{equation*}
\begin{aligned}
A = \begin{bmatrix}
-20 & 22 & 0\\ 0 & 15 & 0 \\  0 &  0  & -3
\end{bmatrix},
\hspace{10pt}B = \mathbb{I}_{3 \times 3}, \hspace{10pt}C = \mathbb{I}_{3 \times 3},
\end{aligned}
\end{equation*}
$f_j(x_{j}) \in \mathbb{R}^{3 \times 1}$ is the system nonlinear vector, $D_j\left(t\right) \in \mathbb{R}^{3 \times 1}$ is the system disturbance vector, $\theta_j \in \mathbb{R}^{3 \times 1}$. Each of $f_j(x_{j}), D_j\left(t\right), \theta_j$  are assumed to be completely unknown.
\begin{equation*}
\begin{aligned}
f_j(x_{j}) = \begin{bmatrix}
a_{1,j}x_{3,j}x_{1,j}+0.2sin(x_{1,j}a_{1,j})\\ -a_{2,j}x_{1,j}x_{3,j}-0.2a_{2,j}cos(a_{2,j}x_{3,j}t)x_{1,j}\\     a_{3,j}x_{1,j}x_{2,j}
\end{bmatrix},
\end{aligned}
\end{equation*}

\begin{equation*}
\begin{aligned}
D_j\left(t\right) = \begin{bmatrix}
1+b_{1,j}sin(b_{1,j}t)\\
1.2cos(b_{2,j}t)\\
sin(0.5b_{3,j}t)+ cos(b_{3,j}t) - 1
\end{bmatrix},
\end{aligned}
\end{equation*}
\begin{equation*}
\begin{aligned}
\theta_j = \begin{bmatrix}
\theta_j^{1} & \theta_j^{2} & \theta_j^{3}
\end{bmatrix},
\end{aligned}
\end{equation*}
\begin{equation*}
\begin{aligned}
\theta_j^{1} &= \begin{bmatrix}
3c_{1,j}sin(0.5t)\\
0.9sin(0.2c_{2,j}t)\\
0.5sin(0.13c_{3,j}t) 
\end{bmatrix},\\
\theta_j^{2} &= \begin{bmatrix}
2c_{1,j}sin(0.4c_{1,j}t)cos(0.3t)\\
2.5sin(0.3c_{2,j}t)+0.3cos\left(t\right)\\
0.6c_{3,j}cos(0.15t)
\end{bmatrix},\\
\theta_j^{3} &= \begin{bmatrix}
0.7sin(0.2c_{1,j}t)\\
1.0sin(0.1c_{2,j}t)\\
1.5cos(0.7c_{3,j}t)+1.6c_{3,j}sin(0.3t)
\end{bmatrix},
\end{aligned}
\end{equation*}
$a, b, c$ are matrices that were selected with different input elements to introduce heterogeneity into the system and therefore different control efforts have to be implemented.
\begin{equation*}
\begin{aligned}
a = \begin{bmatrix}
1.5 & 0.5 & 0.7 & 1.3 & 0.7\\
0.5 & 1.4 & 0.1 & 1.3 & 2.4\\
2.8 & 1.4 & 0.6 & 0.7 & 0.6
\end{bmatrix},\\
b = \begin{bmatrix}
0.5 & 1.5 & 1.1 & 1.6 & 0.3\\
0.7 & 1.2 & 1.3 & 0.5 & 0.3\\
1.1 & 1.4 & 1.6 & 0.6 & 1.0
\end{bmatrix},\\
c = \begin{bmatrix}
1.5 & 2.5 & 0.5 & 1.7 & 0.7\\
0.5 & 1.7 & 1.1 & 0.3 & 0.4\\
0.8 & 0.4 & 2.2 & 0.9 & 1.4
\end{bmatrix},
\end{aligned}
\end{equation*}
The leader's dynamics is selected such that $x_{0} = [3cos(0.7t),
2cos(0.8t), 1.5cos\left(t\right)]^{\top}$. The other parameters of the problem are defined as $\rho_{\infty} = 0.05\times{\bf \underline{1}_{3\times 5}}$, $\rho_{0} = 7\times{\bf \underline{1}_{3\times 5}}$, $l = 7\times{\bf \underline{1}_{3\times 5}}$, $\Gamma = 150\mathbb{I}_{5\times 5}$, $\bar{\delta} = 7$, $\underline{\delta} = 1$, $c = 100$, $k = 0.8$, $\alpha = 50$. Initial conditions of $x(0) = [ 1.6399, 1.6639, -2.1864, 0.1160, -2.7805, -2.2175, -0.1489, $\\
$ 2.2989, -1.3038, 0.5571, -0.5959, 1.6760, -2.4743, 0.0488,$ \\
$0.8288]$. The robustness of the proposed controller against time variant uncertainties in parameters, time-variant disturbances and high nonlinearities are tested in this example considering the formula in \eqref{eq:eq58b}. Fig. \ref{fig:fig9} shows the output performance of the proposed controller for the MIMO case. The control input in the connected graph is shown in Fig. \ref{fig:fig10}. Errors and transformed errors for the three outputs are depicted in Fig. \ref{fig:fig11}, \ref{fig:fig12} and \ref{fig:fig13}. Fig. \ref{fig:fig14} shows the phase plane plot starting from different initial conditions and the synchronization to the desired trajectory. The results demonstrate the performance of the proposed robust controller.

\section{Conclusion}\label{sec_6}
In this paper, a distributed adaptive tracking control of nonlinear uncertain multi-agent systems with prescribed performance is proposed. Under such controller,  the tracking
error is confined from within a predefined large set to a smaller set according to a given performance.  Agents' dynamics were assumed unknown.  The control law is fully distributed based on the fact that the control  of each agent respects the strongly connected graph's topology and includes only the allowed local neighborhood information. The proposed approach guarantees uniform ultimate boundedness for the transformed error.  Simulations include two examples to validate the robustness and smoothness of the proposed controller against highly nonlinear heterogeneous multi-agent system with time-variant uncertain parameters and external disturbances. In future work, control of multi-agents with networks that are  weakly connected or have variable topology will be studied. Under such controllers $L$ and $B$ could time varying and an additional but practical challenges. Systems subjects to actuator failure, saturation or hysteresis will benefit from prescribed performance control and represent interesting areas of further development. Implementation of such control approach on a real system is an area needing research effort.

\section*{Acknowledgment}
The authors would like to acknowledge the support of King Fahd University of Petroleum and Minerals under project Number IN141048.


\bibliographystyle{apacite}
\bibliography{bib_MCC_PPF}

\section*{AUTHOR INFORMATION}
\vspace{10pt}
{\bf Hashim A. Hashim} is a Ph.D. candidate and a Teaching and Research Assistant in Robotics and Control, Department of Electrical and Computer Engineering at the University of Western Ontario, Canada.\\
His current research interests include Stochastic and deterministic filters on SO(3) and SE(3), control of multi-agent systems, control applications and optimization techniques.\\
\underline{Contact Information}: \href{mailto:hmoham33@uwo.ca}{hmoham33@uwo.ca}.
\vspace{40pt}

{\bf Sami El Ferik} is an Associate Professor in Control and Instrumentation, Department of Systems Engineering, at KFUPM. He obtained his B.Sc. in Electrical Engineering from Laval University, Quebec, Canada, and M.Sc. and Ph.D. both in Electrical and Computer Engineering from Polytechnique Montreal, Canada. After the completion of his Ph.D. and Post-doctor positions, he worked with Pratt and Whitney Canada at the Research and Development Center of Systems, Controls, and Accessories. \\
His research interests are in sensing, monitoring, multi-agent systems and nonlinear control with strong multidisciplinary research and applications.\\
\vspace{40pt}

{\bf Frank L. Lewis} National Academy of Inventors, Fellow IEEE, Fellow IFAC, PE Texas, U.K. Chartered Engineer, is a UTA Distinguished Scholar Professor, UTA Distinguished Teaching Professor, and Moncrief-O'Donnell Chair at the University of Texas at Arlington Research Institute. He obtained the Bachelor's Degree in Physics/EE and the MSEE at Rice University, the M.S. in Aeronautical Engineering from the University of West Florida, and the Ph.D. at Ga. Tech. He is author of seven U.S. patents, numerous journal papers, and 14 books.\\

\end{document}